\documentclass{amsart}
\usepackage{amssymb,mathtools}
\usepackage[british]{babel}
\usepackage{enumitem}
\usepackage[latin1]{inputenc}

\usepackage{url}

\newtheorem{theorem}{Theorem}
\numberwithin{theorem}{section}
\newtheorem{corollary}[theorem]{Corollary}
\newtheorem{lemma}[theorem]{Lemma}
\newtheorem{proposition}[theorem]{Proposition}
\theoremstyle{definition}
\newtheorem{definition}[theorem]{Definition}
\newtheorem{remark}[theorem]{Remark}
\newtheorem{example}[theorem]{Example}
\newtheorem{question}[theorem]{Question}

\newtheorem*{claim}{Claim}

\newcommand{\supp}[1]{{\operatorname{supp}_{#1}}}
\newcommand{\rng}{\operatorname{rng}}
\newcommand{\sti}{\Pi^0_2\textsf{-}\mathsf{TI}^\star}
\newcommand{\en}{\operatorname{en}}

\title{Dilators and the reverse mathematics zoo}
\author{Anton Freund}

\address{University of W\"urzburg, Institute of Mathematics, Emil-Fischer-Str.~40, 97074 W\"urz\-burg, Germany}
\email{anton.freund@uni-wuerzburg.de}

\thanks{Funded by the Deutsche Forschungsgemeinschaft (DFG, German Research Foundation) -- Project number 460597863.}

\begin{document}

\begin{abstract}
A predilator is a particularly uniform transformation of linear orders. We have a dilator when the transformation preserves well-foundedness. Over the theory $\mathsf{ACA}_0$ from reverse mathematics, any $\Pi^1_2$-formula is equivalent to the statement that some predilator is a dilator. We show how this completeness result breaks down without arithmetical comprehension: over~$\mathsf{RCA}_0+\mathsf{PA}$, the statements from a large part of the reverse mathematics zoo are not equivalent to some predilator being a dilator.
\end{abstract}

\keywords{Reverse Mathematics, Dilator, Ramsey's Theorem for Pairs}
\subjclass[2020]{03B30, 05D10, 03F15, 03F35}

\maketitle

\section{Introduction}

Consider the transformation of a linear order~$\alpha$ into the order~$2^\alpha$ on terms
\begin{equation*}
2^{\alpha(n-1)}+\ldots+2^{\alpha(0)}\qquad\text{with}\qquad\alpha(0)<\ldots<\alpha(n-1)\text{ in }\alpha,
\end{equation*}
which are compared lexicographically (with $\alpha(n-1)$ taking precedence). Any embedding~$f:\alpha\to\beta$ of linear orders induces an obvious embedding~$2^f:2^\alpha\to 2^\beta$ (where $f$ acts on the exponents), which means that we have a functor on the category of linear orders and embeddings. Let us note that $f\leq g$ implies~$2^f\leq 2^g$ when functions are compared pointwise. Each term $\sigma=2^{\alpha(n-1)}+\ldots+2^{\alpha(0)}\in 2^\alpha$ depends only on a finite subset $\supp{\alpha}(\sigma)=\{\alpha(0),\ldots,\alpha(n-1)\}$ of~$\alpha$. Formally, this corresponds to the fact that all embeddings~$f:\alpha\to\beta$ validate
\begin{equation*}
\left\{\left.\tau\in2^\beta\,\right|\supp{\beta}(\tau)\subseteq\rng(f)\right\}\subseteq\rng\left(2^f\right),
\end{equation*}
where $\rng(g)$ denotes the range or image of a function~$g$. The converse inclusion is entailed by the naturality property
\begin{equation*}
\supp{\beta}\circ 2^f=[f]^{<\omega}\circ\supp{\alpha},
\end{equation*}
where we write $[f]^{<\omega}(a)=\rng(f\restriction a)$ for finite~$a\subseteq\alpha$.

The functorial extension and support functions turn $\alpha\mapsto 2^\alpha$ into a predilator in the sense of J.-Y.~Girard~\cite{girard-pi2}. In Section~\ref{sect:prelim} we will recall the general definition of predilators and some fundamental properties. For the time being, we note that each predilator~$D$ involves a transformation $\alpha\mapsto D(\alpha)$ of linear orders. We only consider countable~$\alpha$ and always assume that this makes~$D(\alpha)$ countable as well. The transformation~$\alpha\mapsto D(\alpha)$ is then computable relative to a certain subset of~$\mathbb N$, which we also denote by~$D$. If this subset is computable, we speak of a computable predilator. When we say that a theory proves a statement about a computable~pre\-dilator, we assume that the latter is given via some program index. Let us note that the coding of predilators by subsets of~$\mathbb N$ makes crucial use of functoriality and the support functions from above, so that the latter are \mbox{well-motivated}.

A predilator~$D$ is a dilator if $D(\alpha)$ is a well-order whenever the same holds for~$\alpha$. In reverse mathematics (see~\cite{simpson09}), the weak base theory~$\mathsf{RCA}_0$ proves that the functorial extension of~$\alpha\mapsto 2^\alpha$ is a predilator. On the other hand, the statement that this predilator is a dilator is equivalent to arithmetical comprehension over~$\mathsf{RCA}_0$, as shown by Girard~\cite{girard87} and J.~Hirst~\cite{hirst94}. The literature now contains many equivalences between natural dilators and important $\Pi^1_2$-principles above~$\mathsf{ACA}_0$, such as infinite iterations of the Turing jump~\cite{marcone-montalban}, arithmetical transfinite recursion \cite{rathjen-weiermann-atr} (originally an unpublished result of H.~Friedman), $\omega$-models of arithmetical transfinite recursion~\cite{rathjen-atr}, $\omega$-models of bar induction~\cite{rathjen-model-bi}, and $\omega$-models of \mbox{$\Pi^1_1$-}comprehension without~\cite{thomson-rathjen-Pi-1-1} and with~\cite{thomson-thesis} bar induction. Since one can quantify over predilators via their codes, one can also consider transformations of higher order, which map dilators to well-orders or again to dilators. This has led to equivalences with the \mbox{$\Pi^1_3$-}principles of $\Pi^1_1$-comprehension~\cite{apw-funct-fast-growing,freund-equivalence,frw-kruskal,uftring-inverse-goodstein} and $\Pi^1_1$-transfinite recursion~\cite{freund-kruskal-friedman,freund-rathjen-iterated-Pi11}.

Girard has shown that the notion of dilator is $\Pi^1_2$-complete. The published proof by D.~Normann (see~\cite[Annex~8.E]{girard-book-part2}) uses the Kleene-Brouwer order and is readily formalized in~$\mathsf{ACA}_0$. This provides some explanation for the results in the previous paragraph, though the equivalences with specific natural dilators require substantial work in each case. Below~$\mathsf{ACA}_0$, the author knows of no proper \mbox{$\Pi^1_2$-}statement that has been characterized by a dilator. The aim of the present paper is to give a systematic explanation for this observation. In particular, we want to understand why research on dilators -- which has found successful applications in the reverse mathematics of better-quasi-orders~\cite{freund-3-bqo} -- has little overlap with investigations into principles from the reverse mathematics zoo (see~\cite{dzhafarov-mummert}), such as Ramsey's theorem for pairs~(denoted $\mathsf{RT}^2_2$; see~\cite{CJS-Ramsey,patey-yokoyama,seetapun-slaman}). To avoid misunderstanding, we stress that proof-theoretic methods and the analysis of well-orders have played a crucial role in the analysis of~$\mathsf{RT}^2_2$ (see~\cite{kolodziejczyk-yokoyama} and the many references discussed there).

We now explain our results in some detail. The following principle will play a central role.

\begin{definition}\label{def:slow-TI}
By slow transfinite $\Pi^0_2$-induction ($\sti$) we mean the statement that $\Pi^0_2$-induction along~$\alpha$ holds for any linear order~$\alpha$ such that $2^\alpha$ is well-founded.
\end{definition}

For comparison, the usual principle of transfinite $\Pi^0_2$-induction asserts that we have $\Pi^0_2$-induction along~$\alpha$ whenever the latter is a well order. It is equivalent to arithmetical comprehension over~$\mathsf{RCA}_0$, as shown by Hirst~\cite{hirst-TI}. Let us agree on the following notation for entailment on a cone of $\omega$-models.

\begin{definition}
In the context of second-order arithmetic, we write
\begin{equation*}
\mathsf T\vDash_\omega\varphi(X_1,\ldots,X_n)
\end{equation*}
if there is a $Y\subseteq\mathbb N$ such that $\mathcal M\vDash\varphi(X_1,\ldots,X_n)$ holds for every $\omega$-model~$\mathcal M\vDash\mathsf T$ with $X_1,\ldots,X_n,Y\in\mathcal M$ (where the set parameters of~$\varphi$ are precisely $X_1,\ldots,X_n$).
\end{definition}

Many unprovability arguments via computability theory yield $\omega$-countermodels. We will later use the following.

\begin{proposition}[C.~Jockusch~\cite{jockusch-ramsey}]\label{prop:jockusch}
We have $\mathsf{RCA}_0\nvDash_\omega\mathsf{RT}^2_2$.
\end{proposition}

The next proposition is established in Section~\ref{sect:sti}. Crucial ideas in the proof are due to F.~Pakhomov and J.~Aguilera (cf.~the paragraph after Lemma~\ref{lem:sti-Pi2-ind}).

\begin{proposition}\label{prop:wkl-sti}
We have~$\mathsf{WKL}_0\nvDash_\omega\sti$.
\end{proposition}

The following dichotomy principle for dilators, which will be proved in Section~\ref{sect:dichot-dil}, is at the center of our approach. We note that the last part of the claim (that the disjunction is exclusive) follows from the previous proposition. This relies on the fact that $\mathsf{RCA}_0\vDash_\omega\varphi$ and $\mathsf{RCA}_0\vDash_\omega\varphi\to\psi$ imply~$\mathsf{RCA}_0\vDash_\omega\psi$, since any two cones of $\omega$-models intersect on a cone.

\begin{theorem}\label{thm:dichot}
For any dilator~$D$, we have
\begin{equation*}
\mathsf{RCA}_0\vDash_\omega\text{``$D$ is a dilator"}\qquad\text{or}\qquad\mathsf{RCA}_0\vDash_\omega\text{``$D$ is a dilator"}\to\sti
\end{equation*}
but not both.
\end{theorem}

The following can also be seen as a conservativity result over $\omega$-models (cf.~the classical result of L.~Harrington, e.\,g., in Corollary~IX.2.6 of~\cite{simpson09}), but we feel that the given formulation as a dichotomy result is more faithful to the constructive content of our argument. An entirely different proof of the corollary (which does not involve Proposition~\ref{prop:wkl-sti} or Theorem~\ref{thm:dichot}) had previously been given by P.~Uftring.

\begin{corollary}\label{cor:wkl-cons}
For any dilator~$D$, we have
\begin{equation*}
\mathsf{RCA}_0\vDash_\omega\text{``$D$ is a dilator"}\qquad\text{or}\qquad\mathsf{WKL}_0\nvDash_\omega\text{``$D$ is a dilator"}.
\end{equation*}
\end{corollary}
\begin{proof}
If the disjunction was false, we could invoke the previous theorem to conclude that we have $\mathsf{WKL}_0\vDash_\omega\sti$, against Proposition~\ref{prop:wkl-sti}.
\end{proof}

So far, we have not been able to strengthen Proposition~\ref{prop:wkl-sti} by including~$\mathsf{RT}^2_2$. Nevertheless, a somewhat indirect approach will yield the following (see Section~\ref{sect:sti}), which can be seen as the main result of our paper.

\begin{theorem}\label{thm:main}
Consider any sentence~$\varphi$ of second-order arithmetic such that we have $\mathsf{WKL}_0+\mathsf{RT}^2_2+\mathsf{PA}\vdash\varphi$ and $\mathsf{RCA}_0\nvDash_\omega\varphi$. For any computable dilator~$D$, we get
\begin{equation*}
\mathsf{RCA}_0+\mathsf{PA}\nvdash\varphi\leftrightarrow\text{``$D$ is a dilator"}.
\end{equation*}
\end{theorem}

We point out that the theorem becomes somewhat trivial when one omits~$\mathsf{PA}$, since all but the most simple dilators rely on~$\Sigma^0_2$-induction (cf.~Proposition~2.2 of~\cite{freund-uftring} and Lemma~\ref{lem:Pi02-ind} below). The theorem clearly applies to a wide range of principles from the reverse mathematics zoo. We only state the following instance.

\begin{corollary}\label{cor:RT22-dil}
For any computable dilator~$D$, we have
\begin{equation*}
\mathsf{RCA}_0+\mathsf{PA}\nvdash\mathsf{RT}^2_2\leftrightarrow\text{``$D$ is a dilator"}.
\end{equation*}
\end{corollary}
\begin{proof}
In view of Proposition~\ref{prop:jockusch} (due to Jockusch), we can apply the previous theorem with $\mathsf{RT}^2_2$ at the place of~$\varphi$.
\end{proof}

Let us mention a result of L.~Ko{\l}odziejczyk and K.~Yokoyama, which says that $\mathsf{WKL}_0+\mathsf{RT}^2_2$ is conservative over~$\mathsf{RCA}_0$ for certain statements about primitive recursive transformations of well-orders below~$\omega^\omega$ (see~\cite[Corollary~3.4]{kolodziejczyk-yokoyama}). This result is not quite comparable with ours, which concerns transformations of arbitrary ordinals rather than standard notation systems. The proof by Ko{\l}odziejczyk and Yokoyama involves new combinatorial results about~$\mathsf{RT}^2_2$, while we combine known results about Ramsey's theorem with a new dichotomy result for dilators.

It is natural to ask whether our Theorem~\ref{thm:main} makes essential use of the specific uniformity properties of dilators or whether it follows from more general results about the reverse mathematics zoo. Being a dilator is a $\forall^1\exists^1\Pi^0_2$-property (where $\forall^1$ and $\exists^1$ refer to set quantifiers). It is known that no statement of this complexity can be equivalent to the cohesive principle, due to a conservation result of D.~Hirschfeldt and R.~Shore (see~\cite[Corollary~2.21]{hirschfeldt-shore}). However, there is still a wide range of principles from the zoo (such as~$\mathsf{RT}^2_2$) for which Theorem~\ref{thm:main} does not follow from results about quantifier complexity. On a somewhat different note, one may wonder how our results depend on the fact that dilators act on linear rather than partial orders (recall that the $\Pi^1_2$-completeness proof for dilators involves the Kleene-Brouwer order). We have no definite answer, though at least the dichotomy in Theorem~\ref{thm:dichot} seems to rely on the specific uniformity properties of dilators and not on linearity alone. Our paper gives rise to several open questions, which are presented at the end of Sections~\ref{sect:dichot-dil} and~\ref{sect:sti}.

\subsection*{Acknowledgements} The author is very grateful for generous and important input from colleagues. Patrick Uftring had previously proved Corollary~\ref{cor:wkl-cons} by entirely different methods. He asked if there are any non-trivial dila\-tors~between $\mathsf{RCA}_0$ and~$\mathsf{ACA}_0$, which inspired this paper. The crucial ideas for the proof of Proposition~\ref{prop:wkl-sti} in Section~\ref{sect:sti} are due to Fedor Pakhomov and Juan Aguilera. A~hint by Lev Beklemishev was important for the proof of Proposition~\ref{prop:induction}. Leszek Ko{\l}odziejczyk and Keita Yokoyama provided information and advice on Ramsey's theorem for pairs and the reverse mathematics zoo.

\section{Preliminaries}\label{sect:prelim}

In this section, we recall the definition of and some fundamental results about dilators. While all mathematical content is due to Girard (see~\cite{girard-pi2} and the introductory~\cite{girard-intro}), we have adapted the presentation to our current aims.

Whenever we speak of a linear order, we assume that its underlying set is a subset of~$\mathbb N$, i.\,e., that we have a representation in second-order arithmetic. Let $\mathsf{LO}$ be the category with the linear orders as objects and the embeddings as morphisms. We define $\mathsf{Nat}$ as the full subcategory with objects $n=\{0,\ldots,n-1\}$ for~$n\in\mathbb N$ (ordered as usual). For order embeddings~$f,g:X\to Y$, we write $f\leq g$ to express that~$f(x)\leq g(x)$ holds for all~$x\in X$. A functor~$D:\mathsf{Nat}\to\mathsf{LO}$ is called monotone if $f\leq g$ implies~$D(f)\leq D(g)$. By $[\cdot]^{<\omega}$ we denote the finite powerset functor on the category of sets, i.\,e., we have $[X]^{<\omega}=\{a\subseteq X\,|\,a\text{ finite}\}$ for any set~$X$ while $[f]^{<\omega}$ for a function~$f$ is defined as in the introduction. We will omit forgetful functors from orders to their underlying sets, and subsets of orders are tacitly considered with the induced order. An example for the following can be found at the very beginning of the present paper.

\begin{definition}\label{def:predil}
A (coded) predilator consists of a monotone functor~$D:\mathsf{Nat}\to\mathsf{LO}$ and a natural transformation~$\supp{}:D\Rightarrow[\cdot]^{<\omega}$ such that the support condition
\begin{equation*}
\{\tau\in D(n)\,|\,\supp{n}(\tau)\subseteq\rng(f)\}\subseteq\rng(D(f))
\end{equation*}
is satisfied for any morphism~$f:m\to n$ of~$\mathsf{Nat}$.
\end{definition}

As noted in the introduction, the converse of the support condition is automatic by naturality, which says precisely that we have the equality in
\begin{equation*}
\supp{n}\left(D(f)(\sigma)\right)=[f]^{<\omega}\left(\supp{m}(\sigma)\right)\subseteq\rng(f).
\end{equation*}
We have defined predilators with domain $\mathsf{Nat}$ rather than~$\mathsf{LO}$ in order to ensure that $D$ and~$\supp{}$ can be represented by subsets of~$\mathbb N$ (see~\cite[Definition~2.4]{freund-computable} for an explicit representation). Let us now explain how an extension to~$\mathsf{LO}$ can be constructed. The idea is that any linear order arises as the direct limit (i.\,e., the union) of its finite suborders, which are isomophic to objects of~$\mathsf{Nat}$. To make this explicit, we write $\en_a:|a|\to a$ for the increasing enumeration of a finite order~$a$ with $|a|$ elements. Each embedding~$f:a\to b$ of finite orders yields a unique morphism~$|f|:|a|\to|b|$ with ${\en_b}\circ|f|=f\circ\en_a$.

\begin{definition}\label{def:predil-extend}
Let $D$ be a coded predilator (with associated transformation~$\supp{}$). With each linear order~$\alpha$ we associate the set
\begin{equation*}
\overline D(\alpha)=\left\{(\sigma,a)\,\left|\,a\in[\alpha]^{<\omega}\text{ and }\sigma\in D(|a|)\text{ with }\supp{|a|}(\sigma)=|a|\right.\right\}.
\end{equation*}
To define a linear order on~$\overline D(\alpha)$, we stipulate
\begin{equation*}
(\sigma,a)\leq(\tau,b)\quad\Leftrightarrow\quad D(|\iota_a|)(\sigma)\leq D(|\iota_b|)(\tau)
\end{equation*}
for the inclusion maps $\iota_a:a\hookrightarrow a\cup b$ and $\iota_b:b\hookrightarrow a\cup b$. When $f:\alpha\to\beta$ is an order embedding, we define $\overline D(f):\overline D(\alpha)\to\overline D(\beta)$ by~$\overline D(f)(\sigma,a)=(\sigma,[f]^{<\omega}(a))$. We also define functions $\operatorname{\overline{supp}}_\alpha:\overline D(\alpha)\to[\alpha]^{<\omega}$ by setting~$\operatorname{\overline{supp}}_\alpha(\sigma,a)=a$.
\end{definition}

Over~$\mathsf{RCA}_0$ one can verify that $\overline D$ and $\operatorname{\overline{supp}}$ form a predilator on~$\mathsf{LO}$, i.\,e., that they satisfy Definition~\ref{def:predil} with $\mathsf{LO}$ at the place of~$\mathsf{Nat}$ (see~\cite{freund-computable}). Conversely, any predilator on~$\mathsf{LO}$ is isomorphic to the extension of some coded predilator, by Girard's normal form theorem~\cite{girard-pi2}. For $\Delta^0_1$-definable dilators on~$\mathsf{LO}$, the isomorphism can be constructed in~$\mathsf{RCA}_0$ (see~\cite[Section~2]{freund-pilot}).

\begin{definition}
A (coded) predilator~$D$ is called a dilator if~$\overline D(\alpha)$ is well-founded for every well order~$\alpha$.
\end{definition}

Let us point out that our definition of dilators is equivalent to the one by Girard. Indeed, the existence of (necessarily unique) support functions is equivalent to the preservation of pullbacks and direct limits (see~\cite[Remark~2.2.2]{freund-thesis}). We also note that Girard's terms~$(\sigma;x_0,\ldots,x_{n-1};\alpha)$ correspond to~$(\sigma,\{x_0,\ldots,x_{n-1}\})$ in our notation. Let us mention that monotonicity on morphisms is automatic when~$\overline D(\omega^\omega)$ is well-founded (see~\cite{girard-pi2} or the presentation in~\cite[Lemma~5.3]{frw-kruskal}).

To ensure that the order on~$\overline D(X)$ is antisymmetric, i.\,e., that its elements have unique representations, we need the condition~$\supp{|a|}(\sigma)=|a|$ from Definition~\ref{def:predil-extend}, which also features in the following.

\begin{definition}
The trace of a predilator~$D$ is given by
\begin{equation*}
\mathsf{Tr}(D)=\{(\sigma,n)\,|\,\sigma\in D(n)\text{ and }\supp{n}(\sigma)=n\}.
\end{equation*}
\end{definition}

In the following result, the point is that $(\sigma,n)\in\mathsf{Tr}(D)$ yields~$(\mu_n(\sigma),n)\in\mathsf{Tr}(E)$, as shown by Girard~\cite{girard-pi2} (see~\cite[Lemma~2.19]{freund-rathjen_derivatives} for the metatheory).

\begin{lemma}[$\mathsf{RCA}_0$]\label{lem:nat-transfo-extends}
Given a natural transformation~$\mu:D\Rightarrow E$ between coded~pre\-dilators, we get a natural transformation~$\overline\mu:\overline D\Rightarrow\overline E$ with~$\overline\mu_\alpha(\sigma,a)=(\mu_{|a|}(\sigma),a)$.
\end{lemma}

To complete our basic setup, we note that each coded predilator~$D$ is isomorphic to the restriction of~$\overline D$ to~$\mathsf{Nat}$. Here the inverse isomorphism maps $(\sigma,a)\in\overline D(n)$ to $D(\en_a^n)(\sigma)\in D(n)$, where $\en_a^n:|a|\to n$ is the embedding with range~$a$  (see~\cite{freund-pilot}). Combined with the previous considerations, we see that predilators on~$\mathsf{LO}$ and coded predilators may be identified on an informal level (though the paragraph before the previous lemma shows that some care is occasionally required).

In the rest of this section, we present some more specialized results of Girard, which will play a central role in the present paper. With respect to the following (originally from~\cite{girard-pi2}), we note that the equalities
\begin{equation*}
\supp{2n}\big(D(e_i)(\sigma)\big)=[e_i]^{<\omega}\big(\supp{n}(\sigma)\big)=\rng(e_i)
\end{equation*}
ensure that the map $i\mapsto D(e_i)(\sigma)$ is injective.

\begin{definition}\label{def:permutation}
Consider a predilator~$D$ and an element~$(\sigma,n)$ of its trace. Let~the embeddings $e_i:n\to 2n$ for $i<n$ be given by~$e_i(i)=2i+1$ and~$e_i(j)=2j$ for~$j\neq i$. We define a permutation~$\pi=\pi_\sigma:n\to n$ by stipulating
\begin{equation*}
D\left(e_{\pi(0)}\right)(\sigma)>\ldots>D\left(e_{\pi(n-1)}\right)(\sigma).
\end{equation*}
When~$a$ is a finite linear order with $n$ elements, we abbreviate~$a^\sigma_i=\en_a(\pi_\sigma(i))$.
\end{definition}

The idea is that an element~$(\sigma,a)\in D(X)$ can be seen as a term with constructor symbol~$\sigma$ and arguments~$\en_a(i)$ for~$i<|a|$ (see the paragraph before Definition~\ref{def:predil-extend}). In view of the following, the permutation~$\pi_\sigma$ determines an order of priority on the arguments, where the importance of position~$\pi_\sigma(i)$ decreases as~$i$ grows. In contrast to Girard, we formulate the result for predilators (rather than just dilators) and with a specific metatheory (cf.~Theorem~3.2.4 of~\cite{girard-pi2}). The proof remains the same despite these changes, but we present it in our notation, not least to demonstrate the specific uniformity properties of predilators that are central to our approach.

\begin{proposition}[$\mathsf{RCA}_0$]\label{prop:permutation-single}
Assume that we have $(\sigma,n)\in\mathsf{Tr}(D)$ for a predilator~$D$. When~$\alpha$ is any linear order, the order on~$\overline D(\alpha)$ satisfies
\begin{equation*}
(\sigma,a)<(\sigma,b)\quad\Leftrightarrow\quad\text{there is a $j<n$ with $a^\sigma_j<b^\sigma_j$ and $a^\sigma_i=b^\sigma_i$ for all~$i<j$}.
\end{equation*}
\end{proposition}
\begin{proof}
Consider the order~$\omega\cdot\alpha$ with elements~$\omega\cdot x+n$ for~$x\in\alpha$ and $n\in\mathbb N$, which are compared according to
\begin{equation*}
\omega\cdot x+n<\omega\cdot x'+n'\quad\Leftrightarrow\quad\text{either $x=x'$ and $n<n'$ or~$x<x'$}.
\end{equation*}
If we have $x<x'$, we write $\omega\cdot x+n\ll\omega\cdot x'+n'$. Note that we have a successor function on~$\omega\cdot\alpha$, which is given by $(\omega\cdot x+n)+1=\omega\cdot x+(n+1)$. Clearly, $\rho\ll\tau$ entails~$\rho+1\ll\tau$. Consider the embedding~$f:\alpha\to\omega\cdot\alpha$ with $f(x)=\omega\cdot x+0$. The equivalence from the proposition retains its truth value when we replace $\alpha$ by~$\omega\cdot\alpha$ and change~$a$ and~$b$ into $[f]^{<\omega}(a)$ and~$[f]^{<\omega}(b)$, respectively. For notational convenience, we assume~$\alpha$ itself is of the form~$\omega\cdot\beta$. Write~$c_i=\en_c(i)$ and $\pi=\pi_\sigma$. For~$j$ as in the proposition, the choice of~$f$ allows us to assume
\begin{equation*}
a^\sigma_j=a_{\pi(j)}\ll\min(a_{\pi(j)+1},b^\sigma_j),
\end{equation*}
where the minimum is evaluated to~$b^\sigma_j$ in case we have~$\pi(j)=n-1$. Under this assumption, we prove the proposition by induction on the number of indices~\mbox{$k<n$} with $a_k>b_k$. Note that it suffices to show the backward implication of our equivalence, as the order on~$\overline D(\alpha)$ is linear. In the base case of the induction, we have~$\en_a\leq\en_b$. To get $(\sigma,a)<(\sigma,b)$, it suffices to recall that predilators are monotone on morphisms. In the induction step, we consider the minimal~$k$ with~$a_k>b_k$. Write $k=\pi(l)$ and note that we get~$j<l$. Let $\overline a$ be the set that results from~$a$ when we replace~$a_j^\sigma$ by~$a^\sigma_j+1$ and $a_k=a^\sigma_l$~by~$b_k$, respectively. In view of $a_{k-1}\leq b_{k-1}<b_k<a_k$ and $a^\sigma_j\ll\min(a_{\pi(j)+1},b_{\pi(j)+1})$, we get $|\overline a|=|a|$~and
\begin{equation*}
\overline a^\sigma_i=\begin{cases}
a^\sigma_j+1 & \text{if $i=j$},\\
b^\sigma_l & \text{if $i=l$},\\
a^\sigma_i & \text{otherwise}.
\end{cases}
\end{equation*}
We thus have $\overline a^\sigma_j<b^\sigma_j$ and $\overline a^\sigma_i=b^\sigma_i$ for~$i<j$, as well as $\overline a^\sigma_j\ll\min(\overline a_{\pi(j)+1},b^\sigma_j)$. Due to the induction hypothesis, it is thus enough to prove the first inequality in
\begin{equation*}
(\sigma,a)<(\sigma,\overline a)<(\sigma,b).
\end{equation*}
We will show that the open inequality follows from~$a^\sigma_j<\overline a^\sigma_j$ and $a^\sigma_l>\overline a^\sigma_l$ with~$j<l$ when we have $a^\sigma_i=\overline a^\sigma_i$ for all~$i\notin\{j,l\}$. Much as above, we may assume that $x<y$ implies~$x\ll y$ for~$x,y\in a\cup\overline a$. We then get an embedding~$g:2n\to\alpha$ by setting
\begin{equation*}
g(2\pi(i))=\begin{cases}
a^\sigma_j & \text{for $i=j$},\\
\overline a^\sigma_l & \text{for $i=l$},\\
a^\sigma_i=\overline a^\sigma_i & \text{for $i\notin\{j,l\}$},
\end{cases}\quad
g(2\pi(i)+1)=\begin{cases}
\overline a^\sigma_j & \text{for $i=j$},\\
a^\sigma_l & \text{for $i=l$},\\
a^\sigma_i+1 & \text{for $i\notin\{j,l\}$}.
\end{cases}
\end{equation*}
For the embeddings~$e_i:n\to 2n$ from Definition~\ref{def:permutation}, we obtain
\begin{equation*}
\left[g\circ e_{\pi(j)}\right]^{<\omega}(n)=\overline a,\qquad\left[g\circ e_{\pi(j)}\right]^{<\omega}(n)=a.
\end{equation*}
Given that $\overline D(g)$ is an embedding, the open claim is thus equivalent to
\begin{equation*}
\left(\sigma,\left[e_{\pi(l)}\right]^{<\omega}(n)\right)<\left(\sigma,\left[e_{\pi(j)}\right]^{<\omega}(n)\right).
\end{equation*}
In view of Definition~\ref{def:predil-extend} (also consider the isomorphism $\overline D(2n)\cong D(2n)$ from the paragraph after Lemma~\ref{lem:nat-transfo-extends}), this inequality amounts to~$D(e_{\pi(l)})(\sigma)<D(e_{\pi(j)})(\sigma)$. The latter holds because of~$j<l$, by the characterization of~$\pi$ in Definition~\ref{def:permutation}.
\end{proof}

In the rest of this section, we deal with comparisons between terms $(\sigma,a)$ and~$(\tau,b)$ that rely on different trace elements.

\begin{definition}\label{def:secure}
Consider a predilator~$D$ as well as two distinct elements $(\sigma,m)$ and $(\tau,n)$ of its trace. We define $P_\sigma^\tau$ as the largest $P\leq\min(m,n)$ with
\begin{equation*}
\pi_\sigma(i)<\pi_\sigma(j)\quad\Leftrightarrow\quad\pi_\tau(i)<\pi_\tau(j)\qquad\text{for all }i,j<P.
\end{equation*}
A number~$p\leq P_\sigma^\tau$ is called secure (relative to $\sigma$ and~$\tau$) if $D(f)(\sigma)<D(g)(\tau)$ has the same truth value for all embeddings $f:m\to m+n$ and~$g:n\to m+n$ such that we have $f(\pi_\sigma(i))=g(\pi_\tau(i))$ for all~$i<p$.
\end{definition}

The following is somewhat stronger than Proposition~6.4.1 of~\cite{girard-pi2} (where the truth value of $D(f)(\sigma)<D(g)(\tau)$ can also depend on the values $f(\pi_\sigma(i))$ for~$i<p$). Except for the first paragraph, our proof follows the one by Girard.

\begin{lemma}[$\mathsf{RCA}_0$]
In the previous definition, $P_\sigma^\tau$ is secure.
\end{lemma}
\begin{proof}
Consider~$f,g$ as in the definition and assume that~$f',g'$ satisfy the same condition. For sufficiently large~$N$, we find embeddings~$h,h':m+n\to N$ such that $h\circ f(\pi_\sigma(i))=h'\circ f'(\pi_\sigma(i))$ holds for all~$i<P=P_\sigma^\tau$. Since~$D(h)$ and $D(h')$ are embeddings, we may replace $f,g$ and~$f',g'$ by $h\circ f,h\circ g$ and $h'\circ f',h'\circ g'$. To simplify notation, we omit $h,h'$ and assume~$f(\pi_\sigma(i))=f'(\pi_\sigma(i))$ for~$i<P$.

If we have $P_\sigma^\tau=m\leq n$, we get $f=g\circ\iota$ with~$\iota(\pi_\sigma(i))=\pi_\tau(i)$, where~$\iota$ is an embedding by the defining property of~$P_\sigma^\tau$. It follows that $D(f)(\sigma)<D(g)(\tau)$ is equivalent to~$D(\iota)(\sigma)=\tau$, which is clearly independent of~$f$ and~$g$. A symmetric argument covers the case where we have $P_\sigma^\tau=n\leq m$. Let us now assume that we have~$P=P_\sigma^\tau<\min(m,n)$. Here the maximality of~$P$ yields a~$j<P$ with
\begin{equation*}
\pi_\sigma(j)<\pi_\sigma(P)\qquad\text{and}\qquad\pi_\tau(j)>\pi_\tau(P)
\end{equation*}
or otherwise $\pi_\sigma(j)>\pi_\sigma(P)$ and $\pi_\tau(j)<\pi_\tau(P)$. Since the two cases are symmetric, we assume that the former applies. Note that we cannot have $D(f)(\sigma)=D(g)(\tau)$, since this would entail
\begin{equation*}
[f]^{<\omega}(m)=[f]^{<\omega}\circ\supp{m}(\sigma)=\supp{m+n}\circ D(f)(\sigma)=[g]^{<\omega}(n),
\end{equation*}
so that we would get~$f=g$ and then~$\sigma=\tau$. By symmetry (exchange $f,g$ with~$f',g'$), we may thus assume that we have~$D(f)(\sigma)<D(g)(\tau)$. The above yields
\begin{equation*}
g(\pi_\tau(P))<g(\pi_\tau(j))=f(\pi_\sigma(j))<f(\pi_\sigma(P)),
\end{equation*}
which remains true when~$f$ and~$g$ are replaced by~$f'$ and~$g'$, respectively. We can now find embeddings~$h,h':m+n\to N$ (not the same as above) such that all $i<P$ validate $h\circ f(\pi_\sigma(i))=h'\circ f(\pi_\sigma(i))$ and we have
\begin{equation*}
h'\circ g(\pi_\tau(P))<h\circ g'(\pi_\tau(P))<h\circ f'(\pi_\sigma(P))<h'\circ f(\pi_\sigma(P)).
\end{equation*}
Note that the middle inequality and the one between the first and the last term are consistent due to the previous observation. To see that, e.\,g., the first inequality is consistent with the choice of values $h\circ f(\pi_\sigma(i))$ for~$i<P$, we note
\begin{equation*}
f(\pi_\sigma(i))<g(\pi_\tau(P))\quad\Leftrightarrow\quad\pi_\tau(i)<\pi_\tau(P)\quad\Leftrightarrow\quad f(\pi_\sigma(i))<g'(\pi_\tau(P)),
\end{equation*}
where the second equivalence relies on the above assumption~$f(\pi_\sigma(i))=f'(\pi_\sigma(i))$. The latter also yields~$h\circ f'(\pi_\sigma(i))=h'\circ f(\pi_\sigma(i))$. Together with the inequality at position~$\pi_\sigma(P)$, we obtain $D(h\circ f')(\sigma)<D(h'\circ f)(\sigma)$ via Proposition~\ref{prop:permutation-single} (observe $a^\sigma_i=h\circ f'$ for~$a=\rng(h\circ f')$ and employ the isomorphism~$D(N)\cong\overline D(N)$ from the paragraph after Lemma~\ref{lem:nat-transfo-extends}). In the same way, we get the last inequality in
\begin{equation*}
D(h\circ f')(\sigma)<D(h'\circ f)(\sigma)<D(h'\circ g)(\tau)<D(h\circ g')(\tau).
\end{equation*}
Since~$D(h)$ is an embedding, we can finally conclude~$D(f')(\sigma)<D(g')(\tau)$.
\end{proof}

Note that the following is well-defined due to the previous lemma.

\begin{definition}\label{def:p-eps}
In the setting of Definition~\ref{def:secure}, we let $p_\sigma^\tau$ be the minimal number that is secure. To choose~$\varepsilon_\sigma^\tau\in\{-1,+1\}$, we declare that we have $\varepsilon_\sigma^\tau={+1}$ when $D(f)(\sigma)<D(g)(\tau)$ holds for some (and hence for all) embeddings~$f:m\to m+n$ and~$g:n\to m+n$ with $f(\pi_\sigma(i))=g(\pi_\tau(i))$ for~$i<p_\sigma^\tau$.
\end{definition}

The following combines Proposition~6.4.2 and Theorem~6.4.5 of~\cite{girard-pi2}.

\begin{theorem}[$\mathsf{RCA}_0$]\label{thm:comp-diff-trace}
Consider a predilator~$D$ and two distinct elements~$(\sigma,m)$ and~$(\tau,n)$ of its trace. For any linear order~$\alpha$ and all finite subsets~$a,b\subseteq\alpha$ with $|a|=m$ and $|b|=n$, the order on~$\overline D(\alpha)$ satisfies
\begin{equation*}
(\sigma,a)<(\tau,b)\quad\Leftrightarrow\quad\begin{cases}
\text{there is~$j<p_\sigma^\tau$ with $a^\sigma_j<b^\tau_j$ and $a^\sigma_i=b^\tau_i$ for all~$i<j$},\\
\text{or we have $a^\sigma_i=b^\tau_i$ for all~$i<p_\sigma^\tau$ and $\varepsilon_\sigma^\tau={+1}$}.
\end{cases}
\end{equation*}
\end{theorem}
\begin{proof}
First note that we have~$p_\tau^\sigma=p_\sigma^\tau$ and~$\varepsilon_\tau^\sigma=-\varepsilon_\sigma^\tau$. So if the right side of the desired equivalence fails, it holds with~$(\sigma,a)$ and~$(\tau,b)$ interchanged. For this reason, it suffices to establish the backward implication. To get $(\sigma,a)<(\tau,b)$, we aim at $D(|\iota_a|)(\sigma)<D(|\iota_b|)(\tau)$ in the notation from Definition~\ref{def:predil-extend}. We note
\begin{equation*}
{\en_{a\cup b}}\circ|\iota_a|(\pi_\sigma(i))=\iota_a\circ\en_a(\pi_\sigma(i))=a^\sigma_i.
\end{equation*}
If the lower line on the right side of our equivalence applies, we can conclude that all~$i<p_\sigma^\tau$ validate~$|\iota_a|(\pi_\sigma(i))=|\iota_b|(\pi_\tau(i))$. Thus $D(|\iota_a|)(\sigma)<D(|\iota_b|)(\tau)$ follows from $\varepsilon_\sigma^\tau={+1}$ (compose $|\iota_a|$ and~$|\iota_b|$ with any embedding of~$|a\cup b|$ into~$m+n$). In the remaining case, we have a $j<p_\sigma^\tau$ with
\begin{equation*}
|\iota_a|(\pi_\sigma(j))<|\iota_b|(\pi_\tau(j))\quad\text{and}\quad |\iota_a|(\pi_\sigma(i))=|\iota_b|(\pi_\tau(i))\text{ for all~$i<j$}.
\end{equation*}
Since $j<p_\sigma^\tau$ is not secure, there are embeddings $f:m\to m+n$ and~$g:n\to m+n$ with $D(f)(\sigma)<D(g)(\tau)$ and~$f(\pi_\sigma(i))=g(\pi_\tau(i))$ for all~$i<j$. When~$N\in\mathbb N$ is sufficiently large, we can find embeddings~$h:|a\cup b|\to N$ and $h':m+n\to N$ such that we have $h\circ|\iota_a|(\pi_\sigma(i))=h'\circ f(\pi_\sigma(i))$ for~$i<j$ as well as
\begin{equation*}
h\circ|\iota_a|(\pi_\sigma(j))<h'\circ f(\pi_\sigma(j))\quad\text{and}\quad h'\circ g(\pi_\tau(j))<h\circ|\iota_b|(\pi_\tau(j)).
\end{equation*}
Let us point out that these two inequalities are consistent with each other since we have $|\iota_a|(\pi_\sigma(j))<|\iota_b|(\pi_\tau(j))$, which is forced in case we have $f(\pi_\sigma(j))<g(\pi_\tau(j))$. For $i<j$ we have
\begin{alignat*}{5}
|\iota_a|(\pi_\sigma(i))&<|\iota_a|(\pi_\sigma(j))\quad&\Leftrightarrow&\quad\pi_\sigma(i)&<\pi_\sigma(j)\quad&\Leftrightarrow&\quad f(\pi_\sigma(i))&<f(\pi_\sigma(j)),\\
f(\pi_\sigma(i))&<g(\pi_\tau(j))\quad&\Leftrightarrow&\quad\pi_\tau(i)&<\pi_\tau(j)\quad&\Leftrightarrow&\quad|\iota_a|(\pi_\sigma(i))&<|\iota_b|(\pi_\tau(j)),
\end{alignat*}
while $\pi_\sigma(i)<\pi_\sigma(j)$ and $\pi_\tau(i)<\pi_\tau(j)$ are equivalent due to~$i,j<p_\sigma^\tau\leq P_\sigma^\tau$. This ensures that all four terms from the two inequalities above fit in the same gap between values $h\circ|\iota_a|(\pi_\sigma(i))=h'\circ f(\pi_\sigma(i))$. As in the previous proof, we can now employ Proposition~\ref{prop:permutation-single} to get
\begin{equation*}
D(h\circ|\iota_a|)(\sigma)<D(h'\circ f)(\sigma)<D(h'\circ g)(\tau)<D(h\circ|\iota_b|)(\tau).
\end{equation*}
This yields~$D(|\iota_a|)(\sigma)<D(|\iota_b|)(\tau)$ and then $(\sigma,a)<(\tau,b)$.
\end{proof}

Let us conclude this section with the following result, which will be needed later (see the proof of Theorem~6.4.6 and Remark~6.4.7 in~\cite{girard-pi2}).

\begin{lemma}[$\mathsf{RCA}_0$]\label{lem:p-inf}
Let $(\rho,k),(\sigma,m)$ and~$(\tau,n)$ be three distinct elements from the trace of a predilator~$D$. We have~$p_\rho^\tau\geq\min(p_\rho^\sigma,p_\sigma^\tau)$.
\end{lemma}
\begin{proof}
We consider an arbitrary~$p<\min(p_\rho^\sigma,p_\sigma^\tau)$ and show that it cannot be secure relative to~$\rho$ and~$\tau$, which forces~$p_\rho^\tau>p$. Note that all~$i,j\leq p$ validate
\begin{equation*}
\pi_\rho(i)<\pi_\rho(j)\quad\Leftrightarrow\quad\pi_\sigma(i)<\pi_\sigma(j)\quad\Leftrightarrow\quad\pi_\tau(i)<\pi_\tau(j).
\end{equation*}
For suitable~$N\in\mathbb N$, we can thus define embeddings~$f,f':k\to N$ and $g,g':n\to N$ as well as~$h:m\to N$ such that
\begin{equation*}
f(\pi_\rho(i))=f'(\pi_\rho(i))=h(\pi_\sigma(i))=g(\pi_\tau(i))=g'(\pi_\tau(i))
\end{equation*}
holds for all~$i<p$ while we have
\begin{equation*}
f(\pi_\rho(p))<h(\pi_\sigma(p))<g(\pi_\tau(p))\quad\text{and}\quad f'(\pi_\rho(p))>h(\pi_\sigma(p))>g'(\pi_\tau(p)).
\end{equation*}
We can now apply the previous theorem (modulo the isomorphism~$D(N)\cong\overline D(N)$ from the paragraph after Lemma~\ref{lem:nat-transfo-extends}) in order to get
\begin{equation*}
D(f)(\rho)<D(h)(\sigma)<D(g)(\tau)\quad\text{and}\quad D(f')(\rho)>D(h)(\sigma)>D(g')(\tau).
\end{equation*}
As promised, this shows that $p$ is not secure for~$\rho$ and~$\tau$.
\end{proof}

\section{A dichotomy result for dilators}\label{sect:dichot-dil}

In the present section, we prove the main part of Theorem~\ref{thm:dichot} from the introduction, namely that at least one of the alternatives from this theorem must hold for any dilator.

The following construction will take centre stage. Let us agree to write~$l(\sigma)$ for the length of a finite sequence~$\sigma$.

\begin{definition}\label{def:D_L}
Consider a linear order~$L=(\mathbb N,\leq_L)$. For each~$n\in\mathbb N$, we define~$\pi^L_n$ as the permutation of~$n=\{0,\ldots,n-1\}$ that is determined by
\begin{equation*}
\pi^L_n(i)\leq_{\mathbb N}\pi^L_n(j)\quad\Leftrightarrow\quad i\leq_L j.
\end{equation*}
For a sequence~$\sigma=\langle\sigma(0),\ldots,\sigma(l(\sigma)-1)\rangle$, we shall abbreviate~$\sigma^L_i=\sigma(\pi^L_{l(\sigma)}(i))$. Given a linear order~$\alpha$, we now set
\begin{equation*}
D_L(\alpha)=\{\langle x_0,\ldots,x_{n-1}\rangle\,|\,x_i\in\alpha\text{ with }x_0<\ldots<x_{n-1}\}.
\end{equation*}
To define a binary relation~$<$ on~$D_L(\alpha)$, we stipulate that all~$\sigma,\tau\in D_L(\alpha)$ validate
\begin{equation*}
\sigma<\tau\quad\Leftrightarrow\quad\begin{cases}
\parbox[t]{\widthof{or we have $l(\sigma)<l(\tau)$ and $\sigma^L_i=\tau^L_i$ for $i<l(\sigma)$.}}{\raggedright either there is $j<\min(l(\sigma),l(\tau))$ with\\[1.2ex] \raggedleft $\sigma^L_j<\tau^L_j$ in~$\alpha$ and $\sigma^L_i=\tau^L_i$ for $i<j$,}\\[5ex]
\text{or we have $l(\sigma)<l(\tau)$ and $\sigma^L_i=\tau^L_i$ for $i<l(\sigma)$.}
\end{cases}
\end{equation*}
For any order embedding~$f:\alpha\to\beta$, we define~$D_L(f):D_L(\alpha)\to D_L(\beta)$ by setting
\begin{equation*}
D_L(f)(\langle x_0,\ldots,x_{n-1}\rangle)=\langle f(x_0),\ldots,f(x_{n-1})\rangle.
\end{equation*}
Finally, we declare that~$\supp{\alpha}:D_L(\alpha)\to[\alpha]^{<\omega}$ is given by
\begin{equation*}
\supp{\alpha}(\langle x_0,\ldots,x_{n-1}\rangle)=\{x_0,\ldots,x_{n-1}\}
\end{equation*}
for any linear order~$\alpha$.
\end{definition}

Let us note that $D_L$ depends not only on the order type of $L$ but also on a given enumeration, i.\,e., on the representation of~$L$ as an order with underlying set~$\mathbb N$. The following is straightforward to verify.

\begin{lemma}[$\mathsf{RCA}_0$]
For any linear order~$L=(\mathbb N,\leq_L)$, the data from the previous definition constitutes a predilator~$D_L$.
\end{lemma}

In the permuted lexicographic order on~$D_L(\alpha)$, the $\pi^L_n(i)$-th smallest argument becomes less important as~$i$ grows. This leads to the following familiar case.

\begin{example}\label{ex:2-dil}
If $\leq_L$ is $\geq_{\mathbb N}$ (so that $L$ has the order type of the negative integers), then~$D_L$ coincides with the dilator~$\alpha\mapsto 2^\alpha$ that was discussed in the introduction.
\end{example}

The next proposition is not needed for any of our main results, but it answers an obvious question. Also, the proof of the forward direction involves an idea that will later reappear in a somewhat more complicated setting.

\begin{proposition}[$\mathsf{ACA}_0$]\label{prop:D_L-ill-founded}
The predilator~$D_L$ is a dilator precisely when the linear order $L=(\mathbb N,\leq_L)$ is ill-founded.
\end{proposition}
\begin{proof}
We begin with the backward direction, which is a variation on the \mbox{well-known} fact that~$\alpha\mapsto 2^\alpha$ preserves well orders. Assume~$\alpha$ is a well order while $\sigma_0,\sigma_1,\ldots$ is an infinitely descending sequence in~$D_L(\alpha)$. We write $\sigma_n=\langle\sigma_n(0),\ldots,\sigma_n(l(n)-1)\rangle$ and put~$\sigma^L_{n,i}=\sigma_n(\pi^L_{l(n)}(i))$ for~$i<l(n)$, where $\pi^L_{l(n)}$ is given as in Definition~\ref{def:D_L}. We inductively find indices~$I(j)$ such that~$n\geq I(j)$ implies $l(n)\geq j$ and~$\sigma^L_{n,i}=\sigma^L_{I(j),i}$ for all~$i<j$. Indeed, if the latter holds, we must have~$l(n)>j$ and~$\sigma^L_{n,j}\geq\sigma^L_{n+1,j}$ for~$n\geq I(j)$, since the~$\sigma_n$ descend. Hence the sequence $I(j)\leq n\mapsto\sigma^L_{n,j}$ in the well order~$\alpha$ must become constant for~$n$ above some~$I(j+1)$. We now show that
\begin{equation*}
\mathbb N\ni j\mapsto\sigma^L_{I(j+1),j}
\end{equation*}
is an embedding of~$L=(\mathbb N,\leq_L)$ into~$\alpha$. Assume that we have~$j\leq_L k$ and that~$I$ is bigger than both $I(j+1)$ and~$I(k+1)$. In view of Definition~\ref{def:D_L}, we get
\begin{equation*}
\sigma^L_{I(j+1),j}=\sigma^L_{I,j}=\sigma_I(\pi^L_{l(I)}(j))\leq\sigma_I(\pi^L_{l(I)}(k))=\sigma^L_{I,k}=\sigma^L_{I(k+1),k}.
\end{equation*}
Given that~$\alpha$ is a well order, it follows that~$L$ cannot be ill-founded.

For the forward direction, we assume that~$L$ is well-founded. We want to show that the infinitely descending sequence
\begin{equation*}
\langle 1\rangle>\langle 0,2\rangle>\langle 0,1,3\rangle>\ldots
\end{equation*}
in the lexicographic order embeds into~$D_L(2\cdot L)$, which will thus be ill-founded  (see the proof of Proposition~\ref{prop:permutation-single} for a definition of product orders). Since~$2\cdot L$ is a well order, this will mean that~$D_L$ cannot be a dilator. We define~$\sigma_n:n+1\to\omega\cdot L$ by
\begin{equation*}
\sigma_n\left(\pi^L_{n+1}(i)\right)=\begin{cases}
2\cdot i+0 & \text{for $i<n$},\\
2\cdot i+1 & \text{for $i=n$}.
\end{cases}
\end{equation*}
Since $\pi^L_{n+1}(i)<\pi^L_{n+1}(j)$ amounts to~$i<_Lj$, the functions~$\sigma_n$ are strictly increasing. We can thus view them as elements~$\langle\sigma_n(0),\ldots,\sigma_n(n)\rangle$ of~$D_L(2\cdot L)$. We have
\begin{equation*}
\sigma_{n+1}\left(\pi^L_{n+2}(n)\right)=\omega\cdot n+0<\omega\cdot n+1=\sigma_n\left(\pi^L_{n+1}(n)\right)
\end{equation*}
as well as $\sigma_{n+1}(\pi^L_{n+2}(i))=\sigma_n(\pi^L_{n+1}(i))$ for all~$i<n$. This yields~$\sigma_{n+1}<\sigma_n$, so that we indeed get a descending sequence in~$D_L$.
\end{proof}

The following definition and result indicate that the dilators~$D_L$ are relevant for a dichotomy result such as Theorem~\ref{thm:dichot}. As preparation for the definition, we note that~$1+\alpha$ has elements~$0$ and~$1+x$ for~$x\in\alpha$, with~$0<1+x<1+y$ for~$x<y$. To turn~$\alpha\mapsto E(\alpha):=\omega\cdot(1+\alpha)$ into a dilator, we set
\begin{alignat*}{3}
E(f)(\omega\cdot 0+n)&=\omega\cdot 0+n,\quad& E(f)(\omega\cdot(1+x)+n)&=\omega\cdot(1+f(x))+n,\\
\operatorname{supp}_\alpha^E(\omega\cdot 0+n)&=\emptyset,\quad& \operatorname{supp}_\alpha^E(\omega\cdot(1+x)+n)&=\{x\}.
\end{alignat*}
When~$D$ is another predilator with support functions~$\operatorname{supp^D_\alpha}$, the usual composition of functors yields a predilator~$D\circ E$ with $\supp{\alpha}:D\circ E(\alpha)\to[\alpha]^{<\omega}$ given by
\begin{equation*}
\supp{\alpha}(\sigma)=\bigcup\left\{\operatorname{supp}^E_\alpha(\rho)\,\left|\,\rho\in\operatorname{supp}^D_{E(\alpha)}(\sigma)\right.\right\}.
\end{equation*}
More precisely, this explains the composition of predilators on~$\mathsf{LO}$ (cf.~the paragraph after Definition~\ref{def:predil-extend}). The case of coded predilators requires some more work, since functors from~$\mathsf{Nat}$ to~$\mathsf{LO}$ cannot be composed directly (see Section~2 of~\cite{freund-rathjen_derivatives}).

\begin{definition}
A linear order~$L=(\mathbb N,\leq_L)$ is called a thread of a predilator~$D$ if there is a natural family of embeddings~$D_L(\alpha)\to D(\omega\cdot(1+\alpha))$.
\end{definition}

In order to explain the occurrence of $\omega\cdot(1+\alpha)$ in the definition, we note that the support sets for~$D_L$ can have any finite cardinality. If we had a natural family of embeddings~$D_L(\alpha)\to D(\alpha)$, the same would need to hold for~$D$, which would be quite restrictive. The choice of $\omega\cdot(1+\alpha)$ is motivated by the following proof.

\begin{theorem}\label{thm:dil-dichotomy}
For any predilator~$D$ such that~$\overline D(\omega^2)$ is well-founded, at least one of the following statements must hold:
\begin{enumerate}[label=\textup{(\roman*)}]
\item We have~$\mathsf{RCA}_0\vDash_\omega\text{``$D$ is a dilator"}$.
\item The predilator~$D$ has a thread.
\end{enumerate}
\end{theorem}
In Corollary~\ref{cor:exclusive}, we will see that the alternative in the theorem is exlusive, i.\,e., that precisely one of statements~(i) and~(ii) holds for each eligible predilator~$D$.
\begin{proof}
Assuming that~(i) fails, we consider an $\omega$-model~$\mathcal M\vDash\mathsf{RCA}_0$ with $D\in\mathcal M$ and a linear order~$\alpha\in\mathcal M$ such that $\mathcal M$ contains no infinitely descending sequence in~$\alpha$ but does contain a sequence
\begin{equation*}
(\sigma(0),a(0))>(\sigma(1),a(1))>\ldots
\end{equation*}
in the order~$\overline D(\alpha)$. Let us abbreviate~$a_{iu}:=a(i)^{\sigma(i)}_u$ for~$u<|a(i)|$ (cf.~Definition~\ref{def:permutation}).

We claim that no trace element~$(\sigma(i),|a(i)|)$ can occur infinitely often. Otherwise, we could assume that some~$(\sigma,n)\in\mathsf{Tr}(D)$ equals~$(\sigma(i),|a(i)|)$ for all~$i\in\mathbb N$. In view of Proposition~\ref{prop:permutation-single}, we would then have an infinitely descending sequence
\begin{equation*}
i\mapsto\alpha^{n-1}\cdot a_{0i}+\ldots+\alpha^0\cdot a_{0,n-1}\in\alpha^n.
\end{equation*}
This would contradict the fact that $\alpha^n$ is well-founded according to~$\mathcal M$, which can be established by $\Pi^0_2$-induction on~$n$ (see, e.\,g., Proposition~9 of~\cite{uftring-etr}). We may now assume that all trace elements~$(\sigma(i),|a(i)|)$ are distinct. This allows us to apply Definition~\ref{def:secure} and the considerations that follow it. Let us agree to write
\begin{equation*}
p_{ij}:=p_{\sigma(i)}^{\sigma(j)}\qquad\text{and}\qquad\varepsilon_{ij}:=\varepsilon_{\sigma(i)}^{\sigma(j)}
\end{equation*}
for~$i<j$ (cf.~Definition~\ref{def:p-eps}).

We now aim to find a sequence of indices~$I(0)<I(1)<\ldots$ (not necessarily in the model~$\mathcal M$) with $p_{I(i),I(j)}\geq i$ for all~$i<j$. Let us inductively assume that we have found a sequence~$I(i)=j(0)<j(1)<\ldots$ in~$\mathcal M$ such that $k<l$ implies~$p_{j(k),j(l)}\geq i$ and additionally $a_{j(k),u}=a_{j(l),u}$ for~$u<i$. Purely for notational convenience, we assume~$j(k)=I(i)+k$. For an inductive construction, it suffices to derive that~$\mathcal M$ contains another sequence~$I(i)<I(i+1):=j(0)<j(1)<\ldots$ such that we have $p_{j(k),j(l)}\geq i+1$ and $a_{j(k),i}=a_{j(l),i}$ for all~$k<l$.

Towards a contradiction, we first assume that all $J$ admit a~$j\geq J$ with $p_{jk}=i$ for all~$k>j$. By passing to a subsequence (which may no longer lie in~$\mathcal M$), we can then assume that all~$j<k$ validate $p_{jk}=i$ as well as $a_{ju}=a_{ku}$ for~$u<i$. In view of Theorem~\ref{thm:comp-diff-trace}, it follows that all inequalities~$\sigma_j>\sigma_k$ in our descending sequence are due to~$\varepsilon_{jk}=-1$ (recall $p_{jk}=p_{kj}$ and~$\varepsilon_{jk}=-\varepsilon_{kj}$ from the proof of the cited theorem). We now find sets $b(j)\in[\omega\cdot(1+i)]^{<\omega}$ with $|b(j)|=|a(j)|$ and
\begin{equation*}
\big\{b(j)^{\sigma(j)}_u\,\big|\,u<i\big\}=\big\{\omega\cdot(1+l)\,\big|\,l<i\big\},
\end{equation*}
which is possible since the gaps below and above the values~$\omega\cdot(1+l)$ accommodate any finite order. In view of Definition~\ref{def:secure}, any~$j<k$ and all $u,w<i=p_{jk}$ satisfy
\begin{equation*}
b(j)^{\sigma(j)}_u<b(j)^{\sigma(j)}_w\quad\Leftrightarrow\quad b(k)^{\sigma(k)}_u<b(k)^{\sigma(k)}_w.
\end{equation*}
Together with the previous equality of sets, this yields
\begin{equation*}
b(j)^{\sigma(j)}_u=b(k)^{\sigma(k)}_u\quad\text{for }i<p_{jk}.
\end{equation*}
As we have~$\varepsilon_{jk}=-1$, we can conclude that $j<k$ implies $(\sigma(j),b(j))>(\sigma(k),b(k))$, due to Theorem~\ref{thm:comp-diff-trace}. This contradicts the assumption that~$\overline D(\omega^2)$ is well-founded.

We now know that any sufficiently large~$j$ admits a~$k>j$ with~$p_{jk}>i$. Hence we can find a sequence~$I(i)<j(0)<j(1)<\ldots$, this time in the model~$\mathcal M$, such that all $k\in\mathbb N$ validate $p_{j(k),j(k+1)}>i$. Due to Lemma~\ref{lem:p-inf}, we even get $p_{j(k),j(l)}>i$ for~$k<l$. As~$(\sigma(j(k)),a(j(k)))$ decreases in~$k$ and $a_{j(k),u}=a_{j(l),u}$ holds for~$u<i$, we can invoke Theorem~\ref{thm:comp-diff-trace} to conclude that $k<l$ entails~$a_{j(k),i}\geq a_{j(l),i}$. Since~$\alpha$ is well-founded according to~$\mathcal M$, we can pass to a tail of our sequence on which we have $a_{j(k),i}=a_{j(l),i}$ for all~$k<l$. This completes the inductive construction. Let us note that the latter also yields~$a_{I(i),u}=a_{I(j),u}$ for~$u<i<j$ (which is relevant for the corollary below but not for the present proof).

If we drop the assumption that the sequence~$i\mapsto(\sigma(i),a(i))$ lies in~$\mathcal M$, we may instead assume that we have~$I(i)=i$, so that $p_{ij}\geq i$ holds for all~$i<j$. Let us recall the permutations~$\pi_{\sigma(i)}$ from Definition~\ref{def:permutation}. We define a linear order~$L=(\mathbb N,\leq_L)$ by stipulating
\begin{equation*}
k\leq_L l\quad\Leftrightarrow\quad\pi_{\sigma(i)}(k)\leq\pi_{\sigma(i)}(l)\qquad\text{for }i>\max\{k,l\},
\end{equation*}
where different~$i$ agree in view of Definition~\ref{def:secure}. For the corollary below, we note that $k\mapsto a_{k+1,k}$ yields an embedding of~$L$ into~$\alpha$. Indeed, we have previously observed that $a_{k+1,k}=a_{Ik}$ holds for~$k<I$. Given~$k<_L l$, we thus get
\begin{equation*}
a_{k+1,k}=a_{Ik}=\en_{a(I)}(\pi_{\sigma(I)}(k))<\en_{a(I)}(\pi_{\sigma(I)}(l))=a_{Il}=a_{l+1,l}
\end{equation*}
with an arbitrary~$I>\max\{k,l\}$. 

In order to establish alternative~(ii) from the state\-ment of the theorem, we first choose embeddings $g_i:|a(i+1)|\to\omega\cdot(1+i)$ with
\begin{equation*}
\{g_i(\pi_{\sigma(i+1)}(k))\,|\,k<i\}=\{\omega\cdot(1+l)\,|\,l<i\},
\end{equation*}
where the values~$g_i(\pi_{\sigma(i+1)}(k))$ for~$i\leq k<|a(i+1)|$ are accommodated in the gaps below and above the elements~$\omega\cdot(1+l)$. An element~$\rho=\langle\rho(0),\ldots,\rho(i-1)\rangle\in D_L(\alpha)$ can be considered as an order embedding~$\rho:i\to\alpha$, which induces an embedding $\overline\rho:\omega\cdot(1+i)\to\omega\cdot(1+\alpha)$ with~$\overline\rho(m)=m$ and $\overline\rho(\omega\cdot(1+l)+m)=\omega\cdot(1+\rho(l))+m$. To satisfy~(ii), we now define~$\eta_\alpha:D_L(\alpha)\to D(\omega\cdot(1+\alpha))$ by setting
\begin{equation*}
\eta_\alpha(\rho)=D(\overline\rho\circ g_i)(\sigma(i+1)).
\end{equation*}
It is straightforward to verify that this definition is natural in the sense that we have $\eta_\beta\circ D_L(f)=D(\omega\cdot(1+f))\circ\eta_\alpha$ for any embedding~$f:\alpha\to\beta$ (reduce this to the fact that we have $\overline{f\circ\rho}=(\omega\cdot(1+f))\cdot\overline\rho$). Let us note that we need only consider finite~$\alpha$ if we work on the level of coded dilators (cf.~Lemma~\ref{lem:nat-transfo-extends}).

To see that~$\eta_\alpha$ is an embedding, recall the permutations~$\pi^L_i$ from Definition~\ref{def:D_L}. The definition of~$L$ ensures that $\pi_i^L(k)\leq\pi_i^L(l)$ and~$\pi_{\sigma(i+1)}(k)\leq\pi_{\sigma(i+1)}(l)$ are~equi\-valent for~$k,l<i$. With $\rho^L_k=\rho(\pi^L_i(k))$ as in Definition~\ref{def:D_L}, we get
\begin{equation*}
\rho^L_k\leq\rho^L_l\quad\Leftrightarrow\quad\overline\rho\circ g_i(\pi_{\sigma(i+1)}(k))\leq\overline\rho\circ g_i(\pi_{\sigma(i+1)}(l))\qquad\text{for }k,l<i.
\end{equation*}
In conjunction with
\begin{align*}
\{\omega\cdot(1+\rho^L_k)\,|\,k<i\}&=\{\omega\cdot(1+\rho(l))\,|\,l<i\}\\
{}&=\{\overline\rho(\omega\cdot(1+l))\,|\,l<i\}=\{\overline\rho\circ g_i(\pi_{\sigma(i+1)}(k))\,|\,k<i\},
\end{align*}
this entails that we have $\omega\cdot(1+\rho^L_k)=\overline\rho\circ g_i(\pi_{\sigma(i+1)}(k))$ for $k<i$.

Now assume we have~$\rho<\tau$ in~$D_L(\alpha)$, for~$\rho$ as before and~$\tau=\langle\tau(0),\ldots\tau(j-1)\rangle$. In view of Definition~\ref{def:D_L}, we first assume that we have~$\rho^L_l<\tau^L_l$ as well as $\rho^L_k=\tau^L_k$ for~$k<l<\min(i,j)$. We get the inequality $\overline\rho\circ g_i(\pi_{\sigma(i+1)}(l))<\overline\rho\circ g_j(\pi_{\sigma(j+1)}(l))$ and equalities below~$l$. By Proposition~\ref{prop:permutation-single} (if $i=j$) or Theorem~\ref{thm:comp-diff-trace} (note that we have $i,j<p_{i+1,j+1}$), we can conclude
\begin{equation*}
\eta_\alpha(\rho)=D(\overline\rho\circ g_i)(\sigma(i+1))<D(\overline\tau\circ g_j)(\sigma(j+1))=\eta_\alpha(\tau).
\end{equation*}
Finally, we consider the case where $\rho<\tau$ holds due to~$\rho^L_k=\tau^L_k$ for all~$k<i<j$. To obtain~$\eta_\alpha(\rho)<\eta_\alpha(\tau)$ via Theorem~\ref{thm:comp-diff-trace}, we need only show
\begin{equation*}
\overline\rho\circ g_i(\pi_{\sigma(i+1)}(i))<\overline\tau\circ g_j(\pi_{\sigma(j+1)}(i))=\omega\cdot(1+\tau^L_i).
\end{equation*}
Write~$g_i(\pi_{\sigma(i+1)}(i))=\omega\cdot(1+l)+m$ with~$l<i$ and possibly~$l=-1$. If the latter holds, the claim is immediate. Now assume~$l\geq 0$ and write $\omega\cdot(1+l)=g_i(\pi_{\sigma(i+1)}(k))$ with~$k<i$. We get $\pi_{\sigma(i+1)}(k)<\pi_{\sigma(i+1)}(i)$ and hence $\pi_{\sigma(j+1)}(k)<\pi_{\sigma(j+1)}(i)$. By
\begin{align*}
\overline\rho\circ g_i(\pi_{\sigma(i+1)}(k))&=\omega\cdot(1+\rho^L_k)=\omega\cdot(1+\tau^L_k)\\
{}&=\overline\tau\circ g_j(\pi_{\sigma(j+1)}(k))<\overline\tau\circ g_j(\pi_{\sigma(j+1)}(i))
\end{align*}
we get the open claim.
\end{proof}

We record the following characterization of dilators, even though it is not needed in the present paper.

\begin{corollary}
A predilator~$D$ is a dilator precisely if~$\overline D(\omega^2)$ is a well order and all threads of~$D$ are ill-founded.
\end{corollary}
\begin{proof}
The forward direction is immediate by Proposition~\ref{prop:D_L-ill-founded}. To establish the back\-ward direction, we assume that~$D$ is no dilator. We then have a well order~$\alpha$ with an infinitely descending sequence~$i\mapsto(\sigma(i),a(i))$ in~$\overline D(\alpha)$. Let us consider an \mbox{$\omega$-}model~$\mathcal M\vDash\mathsf{RCA}_0$ that contains~$D$ and~$\alpha$ as well as our descending sequence. The use of an $\omega$-model is not necessary here, but we prefer to stay close to the previous proof. From the latter, we learn that~$\overline D(\omega^2)$ is ill-founded or that~$D$ has a thread~$L$ that embeds into~$\alpha$ and is thus well-ordered.
\end{proof}

To derive Theorem~\ref{thm:dichot}, we investigate the strength of the statement that~$D_L$ is a dilator for some order~$L$. The following result of Uftring will be used for this purpose. Let us recall that $(1+\alpha)^\gamma$ denotes the order with underlying set
\begin{equation*}
(1+\alpha)^\gamma=\{\langle(x_0,y_0),\ldots,(x_{n-1},y_{n-1})\rangle\,|\,x_0>\ldots>x_{n-1}\text{ in }\gamma\text{ and }y_i\in\alpha\}
\end{equation*}
and lexicographic comparisons (where smaller indices have higher priority and we have $(x_j,y_j)<(x_j',y_j')$ if either $x_j=x_j'$ and $y_j<y_j'$ or~$x_j<x_j'$).

\begin{lemma}[$\mathsf{RCA}_0$; see~\cite{uftring-etr}]\label{lem:Pi02-ind}
For any linear order~$\gamma$, the following are equivalent:
\begin{enumerate}[label=\textup{(\roman*)}]
\item Whenever $\alpha$ is a well order, so is~$(1+\alpha)^\gamma$.
\item The principle of~$\Pi^0_2$-induction along~$\gamma$ holds.
\end{enumerate}
\end{lemma}

The following will yield a connection with slow transfinite induction. Concerning naturality, we note that any order embedding of~$\alpha$ into~$\beta$ yields obvious embeddings of~$(1+\alpha)^\gamma$ into~$(1+\beta)^\gamma$ and of $\alpha\cdot 2^\gamma$ into~$\beta\cdot 2^\gamma$. We do not define these explicitly, because naturality is not needed for the crucial corollary below.

\begin{proposition}[$\mathsf{RCA}_0$]\label{prop:emb-slow-TI}
For arbitrary linear orders~$L=(\mathbb N,\leq_L)$ and~$\gamma$, there is a natural family of embeddings
\begin{equation*}
\eta_\alpha:(1+\alpha)^\gamma\to D_L\left(\alpha\cdot 2^\gamma\right),
\end{equation*}
where~$\alpha$ ranges over all linear orders.
\end{proposition}
\begin{proof}
Our approach is to set
\begin{equation*}
\eta_\alpha(\langle(\gamma_0,\alpha_0),\ldots,(\gamma_{k-1},\alpha_{k-1})\rangle)=\langle e(0),\ldots,e(k-1)\rangle
\end{equation*}
for a suitable embedding~$e:k\to\alpha\cdot 2^\gamma$. For better readability, elements of the order $2^\gamma=(1+\{0\})^\gamma$ are written as $2^{\delta_0}+\ldots+2^{\delta_{n-1}}$ rather than~$\langle(\delta_0,0),\ldots,(\delta_{n-1},0)\rangle$, as in the introduction. Let us recall the permutations~$\pi^L_k$ from Definition~\ref{def:D_L}. To define~$e(\pi^L_k(i))$ by recursion on~$i$, we first set
\begin{equation*}
e(\pi^L_k(0))=\alpha\cdot 2^{\gamma_0}+\alpha_0.
\end{equation*}
Now assume that we have already defined values
\begin{equation*}
e(\pi^L_k(j))=\alpha\cdot\sigma_j+\alpha_j
\end{equation*}
for all~$j\leq i$, where each~$\sigma_j$ is of the form
\begin{equation*}
\sigma_j=2^{\gamma_{j(0)}}+\ldots+2^{\gamma_{j(l)}}\quad\text{with}\quad j(0)<\ldots<j(l)=j.
\end{equation*}
Crucially, there is a unique
\begin{equation*}
\sigma_{i+1}\in\left\{2^{\gamma_{i+1}}\right\}\cup\left\{\left.\sigma_j+2^{\gamma_{i+1}}\,\right|\,j\leq i\right\}
\end{equation*}
such that $e$ remains order preserving when we set
\begin{equation*}
e(\pi^L_k(i+1))=\alpha\cdot\sigma_{i+1}+\alpha_{i+1}.
\end{equation*}
If~$\gamma_j$ is fixed for~$j\leq i$, then $\sigma_{i+1}$ is strictly increasing in~$\gamma_{i+1}$. It is straightforward to conclude that $\eta_\alpha$ is an embedding. Given that the $\sigma_j$ do not depend on the~$\alpha_i$, the construction is clearly natural.
\end{proof}

Let us write~$-\omega$ for the order on the negative integers. Intuitively, the previous proof relies on the fact that $2^{-\omega}$ is a dense linear order, into which any~$L$ can be embedded. While~$-\omega$ does not embed into~$\gamma$, the first components from elements of~$(1+\alpha)^\gamma$ provide descending sequences that allow us to construct partial embeddings in a dynamical way. We note that there is some similarity with the proof of Proposition~\ref{prop:D_L-ill-founded}, though no dynamical approximation was needed in the latter, where we had the (contradictory) assumption that~$L$ is a well order.

\begin{corollary}[$\mathsf{RCA}_0$]\label{cor:D_L-slow-TI}
If $D_L$ is a dilator for some linear order~$L=(\mathbb N,\leq_L)$, then the principle~$\sti$ of slow transfinite $\Pi^0_2$-induction holds (cf.~Definition~\ref{def:slow-TI}).
\end{corollary}
\begin{proof}
In view of Lemma~\ref{lem:Pi02-ind}, it suffices to derive that $\alpha\mapsto(1+\alpha)^\gamma$ preserves well orders whenever~$2^\gamma$ is well-founded. Given that~$D_L$ is a dilator and that~$\alpha$ and~$2^\gamma$ are well orders, the order~$D_L(\alpha\cdot 2^\gamma)$ is well-founded. By the previous proposition, we can conclude that the same holds for~$(1+\alpha)^\gamma$.
\end{proof}

We can now derive the central part of Theorem~\ref{thm:dichot}. The remaining claim, namely that the following disjunction is exclusive, will be establish in the next section.

\begin{theorem}\label{thm:dichot-sti}
For any predilator~$D$ such that~$\overline D(\omega^2)$ is well founded, at least one of the following is true:
\begin{enumerate}[label=\textup{(\roman*)}]
\item We have $\mathsf{RCA}_0\vDash_\omega\text{``$D$ is a dilator"}$.
\item We have $\mathsf{RCA}_0\vDash_\omega\text{``$D$ is a dilator"}\to\sti$.
\end{enumerate}
\end{theorem}
\begin{proof}
Assume that~(i) fails. By Theorem~\ref{thm:dil-dichotomy}, we get a linear order~$L=(\mathbb N,\leq_L)$ and a natural transformation~$\eta$ with components~$\eta_\alpha:D_L(\alpha)\to D(\omega\cdot(1+\alpha))$. When~$\eta$ is considered as a transformation between coded predilators, it is a countable family~of countable order embeddings, which we assume to be coded into a single subset of~$\mathbb N$. In order to establish~(ii), we now consider an $\omega$-model $\mathcal M\vDash\mathsf{RCA}_0+\text{``$D$ is a dilator"}$ that contains $L$ and~$\eta$ (in addition to~$D$). The statement that~$\eta$ is a natural trans\-formation between coded predilators is arithmetical and hence available inside~$\mathcal M$ (cf.~the proof of Theorem~\ref{thm:dil-RT22}). In view of Lemma~\ref{lem:nat-transfo-extends}, we get~$\mathcal M\vDash\text{``$D_L$ is a dilator"}$. The previous corollary allows us to conclude.
\end{proof}

The reader may wonder if Proposition~\ref{prop:emb-slow-TI} and Corollary~\ref{cor:D_L-slow-TI} leave room for improvement. In particular, it is natural to ask whether $2^\alpha$ can be embedded into an order like~$D_L(\omega\cdot (1+\alpha))$ for arbitrary~$L$. By a result that was mentioned in the~introduction, this would allow us to replace~$\sti$ by~$\mathsf{ACA}_0$ in Theorem~\ref{thm:dichot-sti}. We only have the following partial answer. Let us note that the rest of this section does not affect any other part of this paper and may thus be skipped. At the end of the section, we formulate an open question.

\begin{proposition}[$\mathsf{RCA}_0$]\label{prop:scattered}
Assume that $L=(\mathbb N,\leq_L)$ is a scattered linear order of finite Haus\-dorff rank~$r$, i.\,e., that we have an embedding~$h:L\to\mathbb Z^r$. Then there is a natural family of embeddings $\eta_\alpha:2^\alpha\to D_L\big(\alpha\cdot(\omega\cdot(\alpha+1))^r\big)$.
\end{proposition}
\begin{proof}
We intend to set
\begin{equation*}
\eta_\alpha\left(2^{\alpha_0}+\ldots+2^{\alpha_{n-1}}\right)=\langle\sigma(0),\ldots,\sigma(n-1)\rangle
\end{equation*}
for a strictly increasing function~$\sigma:n\to\alpha\cdot(\omega\cdot(\alpha+1))^r$ with values of the form
\begin{equation*}
\sigma\left(\pi^L_n(j)\right)=\alpha\cdot\overline h(j)+\alpha_j.
\end{equation*}
To make~$\eta_\alpha$ an embedding, it is enough to ensure that $\overline h(j)\in(\omega\cdot(\alpha+1))^r$ depends only on~$\alpha_0,\ldots,\alpha_{j-1}$ (i.\,e., not on $\alpha_j,\ldots,\alpha_{n-1}$ and not on~$n$).

The function~$\sigma$ will be strictly increasing if $i<_Lj$ implies~$\overline h(i)<\overline h(j)$. In order to achieve the latter, we rely on the idea that finite descending sequences in~$\alpha$ allow for a dynamic approximation of~$\mathbb Z^r$ by $(\omega\cdot(\alpha+1))^r$ (cf.~the paragraph after the proof of Proposition~\ref{prop:emb-slow-TI}). Write the values of the given embedding~$h:L\to\mathbb Z^r$ as
\begin{equation*}
h(j)=\left\langle q^j_0,\ldots,q^j_{r-1}\right\rangle.
\end{equation*}
For each~$j\in\mathbb N$ we pick an $N_j\geq j$ with $q^j_k\geq -N_j$ for all~$k<r$. Given any~$q\in\mathbb Z$, we define~$j(q)=\min\{j\,|\,q\geq -N_j\}\leq|q|$. Writing~$\alpha_{-1}=\alpha$, we then put
\begin{equation*}
\overline q=\omega\cdot\alpha_{j(q)-1}+N_{j(q)}+q\qquad\text{for } q\in\mathbb Z\text{ with }j(q)\leq n.
\end{equation*}
As we have $j(q^j_k)\leq j$ by construction, we may consider the function
\begin{equation*}
\overline h:n\to(\omega\cdot(\alpha+1))^r\quad\text{with}\quad \overline h(j)=\left\langle \overline{q^j_0},\ldots,\overline{q^j_{r-1}}\right\rangle.
\end{equation*}
Now~$q\mapsto\overline q$ is strictly increasing. To see this, note that~$p<q$ implies~$j(p)\geq j(q)$. If the latter is an equality, we clearly get~$\overline p<\overline q$. If it is a strict inequality, we get the same via $\alpha_{j(p)-1}<\alpha_{j(q)-1}$. Given that~$h$ is an embedding, it follows that~$i<_L j$ implies~$\overline h(i)<\overline h(j)$, as desired. Again due to $j(q^j_k)\leq j$, we can confirm that~$\overline h(j)$ depends on $\alpha_0,\ldots,\alpha_{j-1}$ only. The construction is clearly natural in~$\alpha$.
\end{proof}

Recall that, over~$\mathsf{RCA}_0$, arithmetical comprehension is equivalent to the statement that~$\alpha\mapsto 2^\alpha$ preserves well orders (see~\cite{girard87,hirst94}). Together with the previous propo\-sition, this yields the following (cf.~the proof of Theorem~\ref{thm:dichot-sti}).

\begin{corollary}
For any predilator~$D$ with a thread of finite Hausdorff rank, we have $\mathsf{RCA}_0\vDash\text{``$D$ is a dilator"}\to\mathsf{ACA_0}$.
\end{corollary}

In the following, we explain the restriction to finite Hausdorff ranks.

\begin{remark}
In Proposition~\ref{prop:scattered}, one can replace~$\mathbb Z^r$ by $(1+\mathbb Z)^\gamma$ for an infinite ordinal~$\gamma$. However, this does not yield a result about infinite Hausdorff ranks. The reason lies in different definitions of exponentiation for linear orders. When we write $(1+\mathbb Z)^\gamma$, we refer to the definition in the paragraph before Lemma~\ref{lem:Pi02-ind}. On the other hand, the definition of Hausdorff rank (see~\cite{rosenstein-linear-orders}) takes $\mathbb Z^\gamma$ to be the order on functions~$f:\gamma\to\mathbb Z$ with only finitely many non-zero values, where we have~$f<g$ if $f(\beta)<g(\beta)$ holds for the biggest~$\beta$ on which~$f$ and~$g$ differ. If we only record non-zero values, such functions can also be represented by finite sequences $\langle(\gamma_0,q_0),\ldots,(\gamma_{n-1},q_{n-1})\rangle$, but we have $\langle(0,-1)\rangle>\langle(1,-1)\rangle>\langle(2,-1)\rangle>\ldots$ in the order $\mathbb Z^\omega$. For infinite~$\gamma$, this makes it hard to approximate $\mathbb Z^\gamma$ by well orders, so that our proof of Proposition~\ref{prop:scattered} does not extend to this case.
\end{remark}

To conclude this section, we formulate an open problem.

\begin{question}
How strong is the statement that $D_L$ is a dilator for some linear order~$L=(\mathbb N,\leq_L)$? Is it stronger than~$\sti$ (cf.~Corollary~\ref{cor:D_L-slow-TI})? Is it weaker than arithmetical comprehension? What happens in the extreme case~$L\cong\mathbb Q$?
\end{question}

\section{On slow transfinite induction}\label{sect:sti}

In the present section, we prove bounds on the strength of slow transfinite induction and derive the results that were stated in the introduction.

Let us begin with a straight\-forward observation. As usual, we write $\mathsf{I\Sigma}^0_n$ to refer to $\Sigma^0_n$-induction along the natural numbers.

\begin{lemma}\label{lem:sti-Pi2-ind}
\textup{(a)} We have $\mathsf{RCA}_0+\sti\vdash\mathsf{I\Sigma}^0_2$.
\begin{enumerate}[label=\textup{(\alph*)}]\setcounter{enumi}{1}
\item Over~$\mathsf{RCA}_0+\mathsf{I\Sigma}^0_2$, the following is equivalent to~$\sti$: a linear order~$\gamma$ such that~$2^\gamma$ is well-founded admits no infinitely descending $\Delta^0_2$-sequence.
\end{enumerate}
\end{lemma}
\begin{proof}
(a) The order type~$\omega$ of the natural numbers is isomorphic to~$2^\omega$ (think of binary notation). This fact and the well-foundedness of~$\omega$ is  recognized in~$\mathsf{RCA}_0$. Given~$\sti$, we thus get~$\Pi^0_2$-induction along~$\mathbb N$. It is well-known that the latter is equivalent to~$\mathsf{I\Sigma}^0_2$ (see, e.\,g., Section~I.2 of~\cite{hajek91}).

(b) It suffices to show that $\Pi^0_2$-induction along~$\gamma$ holds precisely if no $\Delta^0_2$-sequence in~$\gamma$ descends infinitely. For the forward direction, we assume that $F:\mathbb N\to\gamma$ is~$\Delta^0_2$ with $F(n)>F(n+1)$ for all~$n\in\mathbb N$. We obtain $F(n)>x$ for all~$n\in\mathbb N$ by a straightforward induction on~$x\in\gamma$, which clearly yields a contradiction. For the backward direction, we assume that induction along~$\gamma$ fails for a $\Pi^0_2$-formula~$\psi(x)$, i.\,e., that we have
\begin{equation*}
\exists x\in\gamma\,\neg\psi(x)\quad\text{and}\quad\forall x\in\gamma(\neg\psi(x)\to\exists y<x\,\neg\psi(y)).
\end{equation*}
The idea is to construct~$x_0>x_1>\ldots$ in~$\gamma$ such that $\neg\psi(x_i)$ holds for all~$i\in\mathbb N$. To achieve this in our metatheory, we write $\psi(x)$ as~$\forall n\,\varphi(x,n)$ where~$\varphi$ is~$\Sigma^0_1$. Let~$S$ be the $\Delta^0_2$-collection of all finite sequences~$\sigma=\langle\sigma_0,\ldots,\sigma_{k-1}\rangle$ such that each~$\sigma_i$ is the minimal code of a pair~$\langle x_i,n_i\rangle$ with~$\neg\varphi(x_i,n_i)$ and either~$i=0$ or~$x_i<x_{i-1}$. Using~$\mathsf{I\Sigma}^0_2$, one readily shows that~$S$ contains a (clearly unique) sequence of each length. In order to obtain a $\Delta^0_2$-function $F:\mathbb N\to\gamma$ with $F(n)<F(n-1)$, we now declare that~$F(n)$ is the first component of the $n$-th pair in some (or any) suitably long sequence from~$S$.
\end{proof}

We now provide the proof for a proposition from the introduction, which asserts that we have $\mathsf{WKL}_0\nvDash_\omega\sti$. It was found at \emph{Trends in Proof Theory}~2024. The crucial ideas for the second and third paragraph of the proof were suggested by F.~Pakhomov and J.~Aguilera, respectively.

\begin{proof}[Proof of Proposition~\ref{prop:wkl-sti}]
Given any~$Z\subseteq\mathbb N$, we need to find an~$\omega$-model~$\mathcal M\vDash\mathsf{WKL}_0$ that contains~$Z$ and falsifies~$\sti$. Let us anticipate that the following arguments relativize, which allows us to assume~$Z=\emptyset$. We find a low set~$X$ and an $\omega$-model $\mathcal M\vDash\mathsf{WKL}_0$ such that all sets in~$\mathcal M$ are computable from~$X$ (take a low $\mathsf{PA}$-degree, e.\,g.~via Theorems~3.17 and~4.22 of~\cite{hirschfeldt-slicing}). In view of the previous lemma, it is enough to find a computable linear order~$\gamma$ such that $2^\gamma$ has no~$X$-computable~descending sequence while $\gamma$ contains a descending sequence that is computable in $0'$.

For an $X$-computable linear order~$\gamma$, it is known that an~$X$-computable descending sequence in~$2^\gamma$ yields an $X'$-computable descending sequence in~$\gamma$ itself (see Theorem~3.1 of~\cite{marcone-montalban}). To formulate an effective version of this fact, we write $\{e\}_X(n)$ and $\{e\}_X'(n)$ for the result (possibly undefined) of the $e$-th programme on input~$n$ with oracle~$X$ and~$X'$, respectively. The proof in the cited reference yields a computable function~$f:\mathbb N^2\to\mathbb N$ with the following properties:
\begin{itemize}
\item If $\{e\}_X$ is the characteristic function of the graph of a linear order~$\gamma$ and $\{s\}_X$ is a strictly descending sequence in~$\gamma$, then $\{f(e,s)\}_X'$ is a strictly~descending sequence in~$\gamma$.
\item The function $\{f(e,s)\}_X'$ is total for all~$e,s\in\mathbb N$.
\end{itemize}
Now the idea is to construct a computable linear order~$\gamma$ with a $0'$-computable~descending sequence such that every sequence~$x_0>x_1>\ldots$ in~$\gamma$ validates
\begin{equation*}
\{f(e,s)\}_X'(n)<_{\mathbb N} x_n\quad\text{for all }n=\langle e,s\rangle\in\mathbb N.
\end{equation*}
If~$e_0$ and $s_0$ were $X$-indices for~$\gamma$ and a descending sequence in~$2^\gamma$, the sequence that is given by $x_n=\{f(e_0,s_0)\}_X'(n)\in\gamma$ would yield a contradiction at~$n=\langle e_0,s_0\rangle$. Thus $2^\gamma$ will have to be $X$-computably well-founded, as demanded above.

Given that~$X$ is low, we get a total $0'$-computable function $D:\mathbb N\to\mathbb N$ by setting
\begin{equation*}
D(n)=\{f(e,s)\}_X'(n)+1\quad\text{for }n=\langle e,s\rangle.
\end{equation*}
By the limit lemma (see Section~3.6 of~\cite{soare-computability}), we get $D(n)=\lim_k d(k,n)$ for a~computable function~$d:\mathbb N^2\to\mathbb N$, where we can assume~$d(k,n)\leq k$. Let us write $K(n)$ for the smallest~$K$ such that $D(n)=d(k,n)$ holds for all~$k\geq K$. We point out that we get~$D(n)\leq K(n)$. Define~$T$ as the tree of sequences~$\langle k(0),\ldots,k(l-1)\rangle$ such that the following holds for~$n<l$:
\begin{enumerate}[label=(\roman*)]
\item If $k(n)$ is non-zero, then we have $d(k(n)-1,n)\neq d(k(n),n)$.
\item We have $d(k(n),n)=d(k,n)$ for~$k(n)<k<l$.
\end{enumerate}
Condition~(i) ensures~$k(n)\leq K(n)$. The latter must be an equality when the given sequence lies on a path of~$T$, by condition~(ii) and the minimality of~$K(n)$. Thus~$T$ has a unique path, which consists of the rightmost sequences. Let~$\gamma$ be the Kleene-Brouwer order on~$T$ without its root. Any descending sequence~$x_0>x_1>\ldots$ in~$\gamma$ must lie entirely within the branch. Assuming that the code of a sequence~ma\-jor\-izes its entries, we thus get $x_n\geq K(n)\geq D(n)$, as needed.
\end{proof}

Let us note that Corollary~\ref{cor:wkl-cons} from the introduction is now established as well. As promised in the previous section, we get the following.

\begin{corollary}\label{cor:exclusive}
The alternatives in Theorems~\ref{thm:dil-dichotomy} and~\ref{thm:dichot-sti} are exclusive, i.\,e., in each of these two results and for each predilator~$D$, statements~(i) and~(ii) cannot hold simultaneously.
\end{corollary}
\begin{proof}
Statement~(i) is the same in both results, and statement~(ii) of Theorem~\ref{thm:dil-dichotomy} implies (ii) of Theorem~\ref{thm:dichot-sti}, by the proof of the latter. Thus it suffices to show that the alternative in the second result is exclusive. If this was not the case, we would get~$\mathsf{RCA}_0\vDash_\omega\sti$, against the result of the previous proof.
\end{proof}

Our next aim is to derive a result that involves Ramsey's theorem for pairs and two colours ($\mathsf{RT}^2_2$). We begin with some preparations. Let us recall that $\varepsilon_0$, which is the proof-theoretic ordinal of Peano arithmetic ($\mathsf{PA}$), has a standard representation as an ordered system of terms that are recursively generated by the clause
\begin{equation*}
t(0)\succeq\ldots\succeq t(n-1)\text{ in }\varepsilon_0\qquad\Rightarrow\qquad\omega^{t(0)}+\ldots+\omega^{t(n-1)}\in\varepsilon_0.
\end{equation*}
This already refers to the order~$\prec$ on~$\varepsilon_0$, which is lexicographic with respect to the recursively defined order on exponents (see, e.\,g., Section~3.3.1 of~\cite{pohlers-proof-theory} for more details). The following proof, which adapts a classical argument due to G.~Kreisel and A.~L\'evy~\cite{kreisel68}, was pointed out to the author by L.~Beklemishev. We write~$\exists^1\Pi^0_n$ for the class of formulas~$\exists X\subseteq\mathbb N\,\varphi$ such that $\varphi$ is~$\Pi^0_n$. In the context of second-order arithmetic, $\mathsf{PA}$ denotes the schema of arithmetical induction with set parameters.

\begin{proposition}\label{prop:induction}
When~$\psi$ is any true $\exists^1\Pi^0_3$-sentence, the theory $\mathsf{RCA}_0+\mathsf{PA}+\psi$ cannot prove that the well-foundedness of~$\varepsilon_0$ entails $\Pi^0_2$-induction along~$\varepsilon_0$.
\end{proposition}
\begin{proof}
We begin with a reduction to a corresponding first-order result. Let us extend the language of first-order arithmetic by a unary predicate symbol~$U$. We will use~$U$ to account for the second-order quantifier in~$\psi$. For a $\Pi^1_1$-sentence~$\forall X\,\varphi$ of second-order arithmetic, we define~$\varphi^U$ as the first-order sentence that results from~$\varphi$ when all (positive and negative) subformulas $t\in X$ are replaced by~$U(t)$. In the following, it is understood that~$U$ may occur in the induction formulas of first-order Peano arithmetic.

We consider transfinite induction only along (our standard representation of) the ordinal~$\varepsilon_0$. When $\Gamma$ is a class of formulas (each with a distinguished induction variable), we write $\Gamma\textsf{-TI}$ for the collection of formulas
\begin{equation*}
\forall\gamma\in\varepsilon_0\big(\forall\beta\prec\gamma\,\varphi(\beta)\to\varphi(\gamma)\big)\to\forall\alpha\in\varepsilon_0\,\varphi(\alpha)
\end{equation*}
with $\varphi\in\Gamma$. We will be particularly interested in the case where~$\Gamma$ is the class of $\Pi_n$-formulas with~$n\in\{1,2\}$. Let us agree that such formulas may contain~$U$. We have the following conservation result.

\begin{claim}
For any $\Pi^1_1$-sentence~$\forall X\,\varphi$, we have
\begin{equation*}
\mathsf{RCA}_0+\mathsf{PA}+\text{``$\varepsilon_0$ is well-founded"}\vdash\forall X\,\varphi\quad\Rightarrow\quad\mathsf{PA}+\Pi_1\text{-}\mathsf{TI}\vdash\varphi^U.
\end{equation*}
\end{claim}

When $\mathsf{PA}$ is replaced by~$\mathsf{I\Sigma}_1$, the result coincides with Theorem~4 of~\cite{rathjen-carlucci-mainardi}. The proof extends the usual model-theoretic argument for the conservativity of~$\mathsf{RCA}_0$ over~$\mathsf{I\Sigma_1}$. Aiming at the contrapositive of our implication, we consider a model~$\mathcal M$ of the theory~$\mathsf{PA}+\Pi_1\text{-}\mathsf{TI}+\neg\varphi^U$. Let~$\mathcal S$ consist of the $\Delta_1$-definable subsets of the universe~$M$ of~$\mathcal M$ (where definitions may involve~$U$). If a subset of~$M$ is \mbox{$\Delta^0_1$-}definable with number and set parameters over~$(\mathcal M,\mathcal S)$, it is $\Delta_1$-definable with number parameters over~$\mathcal M$ (see~\cite[Lemma~IX.1.8]{simpson09}). The analogous claim for first-order definability follows. Thus $(\mathcal M,\mathcal S)$ validates $\mathsf{RCA}_0+\mathsf{PA}$ and clearly also~$\neg\forall X\,\varphi$. As in~\cite{rathjen-carlucci-mainardi}, ill-foundedness of~$\varepsilon_0$ in~$(\mathcal M,\mathcal S)$ would yield a failure of $\Pi_1\text{-}\mathsf{TI}$ in~$\mathcal M$. We will use the claim above in conjunction with the following.

\begin{claim}
For any~$\Sigma_4$-sentence~$\varphi$ that is true under some interpretation of~$U$, we have 
\begin{equation*}
\mathsf{PA}+\varphi+\Pi_1\text{-}\mathsf{TI}\nvdash\Pi_2\text{-}\mathsf{TI}.
\end{equation*}
\end{claim}

By a classical result of Kreisel and L\'evy~\cite{kreisel68}, induction along~$\varepsilon_0$ is equivalent to uniform reflection over~$\mathsf{PA}$. A closer look at the proof shows that $\Pi_2$-induction corresponds to~$\Pi_4$-reflection. To formalize the latter, we use a $\Pi_4$-formula~$\mathsf{Tr}_4$ that defines truth for~$\Pi_4$-sentences (see~\cite[Section~I.1]{hajek91}). We also have a $\Sigma_1$-formula~$\mathsf{Pr}_{\mathsf{PA}}$ that expresses provability in~$\mathsf{PA}$. Our reflection principle can now be stated as
\begin{equation*}
\mathsf{RFN}_4=\forall\theta\,\big(\text{``$\theta$ is a $\Pi_4$-sentence"}\land\mathsf{Pr}_{\mathsf{PA}}(\theta)\to\mathsf{Tr}_4(\theta)\big).
\end{equation*}
The result of Kreisel and L\'evy (refined with respect to formula complexity) says
\begin{equation*}
\mathsf{PA}+\Pi_2\text{-}\mathsf{TI}\vdash\mathsf{RFN}_4.
\end{equation*}
For a detailed proof via ordinal analysis, we refer to Section~5 of~\cite{freund-Sigma2-incompleteness} (also see the alternative approach by Beklemishev~\cite{beklemishev03,beklemishev-provability-algebras}). There it is shown that parameter-free $\Pi_1$-induction implies local~$\Sigma_2$-reflection. No new issues arise when we admit parameters and go from~$\Pi_2$-induction to uniform~$\Sigma_3$-reflection, which is known to be equivalent to uniform~$\Pi_4$ reflection. The addition of the predicate symbol~$U$ amounts to a harmless relativization. We can express $\Pi_1\text{-}\mathsf{TI}$ by a single $\Pi_3$-sentence (again via a partial truth definition). Let $\theta$ be the $\Sigma_4$-sentence $\varphi\land\Pi_1\text{-}\mathsf{TI}$. If the proposition was false, we would get $\mathsf{PA}+\theta\vdash\mathsf{Pr}_{\mathsf{PA}}(\neg\theta)\to\neg\theta$. The assumption on~$\varphi$ ensures that~$\mathsf{PA}+\theta$ is consistent. So we would have a contradiction with G\"odel's second incompleteness theorem. This completes the proof of our second claim.

We now write $\exists X\,\psi_0$ for the $\exists^1\Pi^0_3$-formula~$\psi$ from the statement of the proposition. For any set variable~$Y$, let $\Pi_2\textsf{-TI}_Y$ be the second-order formula that results from $\Pi_2\textsf{-TI}$ when each subformula~$U(t)$ is replaced by~$t\in Y$. Note that $\forall Y\,\Pi_2\textsf{-TI}_Y$ is a consequence of the second-order version of $\Pi^0_2$-induction along~$\varepsilon_0$. Assuming that the proposition is false, we thus get
\begin{equation*}
\mathsf{RCA}_0+\mathsf{PA}+\text{``$\varepsilon_0$ is well-founded"}\vdash\exists X\,\psi_0\to\forall Y\,\Pi_2\textsf{-TI}_Y.
\end{equation*}
Here the conclusion can be weakened to $\forall X(\psi_0\to\Pi_2\textsf{-TI}_X)$. Due to the first claim above, we can infer
\begin{equation*}
\mathsf{PA}+\Pi_1\textsf{-TI}\vdash\psi_0^U\to\Pi_2\text{-TI}.
\end{equation*}
But this contradicts the second claim.
\end{proof}

To get the following, it suffices to combine the previous proposition with a recent result of Q.~Le~Hou\'erou, L.~Patey and K.~Yokoyama~\cite{houerou-patey-yokoyama}.

\begin{corollary}\label{cor:RT22-sTI}
When~$\psi$ is a true~$\exists^1\Pi^0_3$-sentence, we have
\begin{equation*}
\mathsf{WKL}_0+\mathsf{RT}^2_2+\mathsf{PA}+\psi\nvdash\sti.
\end{equation*}
\end{corollary}
\begin{proof}
It is straightforward to show that $2^{\varepsilon_0}$ embeds into~$\varepsilon_0$ (in fact by an isomorphism). Assuming that the corollary fails, we infer that the theory
\begin{equation*}
\mathsf{WKL}_0+\mathsf{RT}^2_2+\mathsf{PA}+\text{``$\varepsilon_0$ is well-founded"}
\end{equation*}
proves that $\psi$ implies~$\Pi^0_2$-induction along~$\varepsilon_0$. This implication is a $\Pi^1_1$-statement. We can thus use Theorem~1.6 of~\cite{houerou-patey-yokoyama} to infer that the same implication is provable in $\mathsf{RCA}_0+\mathsf{PA}+\text{``$\varepsilon_0$ is well-founded"}$. But this contradicts the previous proposition.
\end{proof}

Finally, we derive the following result, which is a reformulation of Theorem~\ref{thm:main} from the introduction.

\begin{theorem}\label{thm:dil-RT22}
For any computable dilator~$D$, we have
\begin{equation*}
\mathsf{RCA}_0\vDash_\omega\text{``$D$ is a dilator"}\qquad\text{or}\qquad\mathsf{WKL}_0+\mathsf{RT}^2_2+\mathsf{PA}\nvdash\text{``$D$ is a dilator"}.
\end{equation*}
\end{theorem}
\begin{proof}
Towards a contradiction, we assume the theorem fails for~$D$. Then we have
\begin{equation*}
\mathsf{WKL}_0+\mathsf{RT}^2_2+\mathsf{PA}\vdash\text{``$D$ is a dilator"}.
\end{equation*}
Also, Theorem~\ref{thm:dil-dichotomy} ensures that $D$ has a thread, i.\,e., that the following statement~$\psi$ is true: there is a natural family of order embeddings~$\eta_n:D_L(n)\to D(\omega\cdot(1+n))$ with~$n\in\mathbb N$ for some linear order~$L=(\mathbb N,\leq_L)$. In the proof of Theorem~\ref{thm:dichot-sti} we have used that~$\psi$ is~$\Sigma^1_1$. Here we note that it is actually~$\exists^1\Pi^0_2$. Indeed, the order~$D_L(n)$ and thus the embedding~$\eta_n$ is finite for each~$n\in\mathbb N$. Our natural family can thus be seen as a function~$n\mapsto\eta_n\in\mathbb N$, and it takes a $\Pi^0_2$-formula to express that this function is total. To say that the family is natural, we need to quantify over all morphisms~$f:m\to n$ with~$m,n\in\mathbb N$, which amounts to a first-order quantification. For each~$f$, the naturality condition~$\eta_n\circ D_L(f)=D(\omega\cdot(1+f))\circ\eta_m$ asserts only finitely many equalities and is thus $\Delta^0_1$. By Lemma~\ref{lem:nat-transfo-extends}, our natural family yields an embedding~$\overline{D_L}(\alpha)\to\overline D(\omega\cdot(1+\alpha))$ for any linear order~$\alpha$, provably in~$\mathsf{RCA}_0$. We thus get
\begin{equation*}
\mathsf{WKL}_0+\mathsf{RT}^2_2+\mathsf{PA}+\psi\vdash\text{``$D_L$ is a dilator for some linear order~$L=(\mathbb N,\leq_L)$}".
\end{equation*}
By Corollary~\ref{cor:D_L-slow-TI}, the same theory proves~$\sti$, against Corollary~\ref{cor:RT22-sTI}.
\end{proof}

The following complements the question from the end of the previous section. Yokoyama has pointed out that a positive answer to~(c) would yield a more canonical version of the previous theorem.

\begin{question}
\textup{(a)} Do we have~$\mathsf{WKL}_0+\mathsf{RT}^2_2\nvDash_\omega\sti$?
\begin{enumerate}[label=\textup{(\alph*)}]\setcounter{enumi}{1}
\item Is $\sti$ strictly weaker than the principle of arithmetical comprehension? Is it $\Pi^1_1$-conservative over~$\mathsf{RCA}_0+\mathsf{I\Sigma}^0_2$ (cf.~Lemma~\ref{lem:sti-Pi2-ind})?
\item Is $\mathsf{WKL}_0+\mathsf{RT}^2_2+\mathsf{PA}$ conservative over~$\mathsf{RCA}_0+\mathsf{PA}$ for statements of the form ``$D$ is a dilator" with computable~$D$?
\end{enumerate}
\end{question}

\bibliographystyle{amsplain}
\bibliography{RM-Zoo-Dilators}

\end{document}